\documentclass[10pt]{amsart}
\usepackage{amsmath,amssymb,latexsym,amsthm, amscd, mathrsfs}

\theoremstyle{definition}
\newtheorem{example}[subsection]{Example}
\theoremstyle{plain}
\newtheorem{prop}[subsection]{Proposition}
\newtheorem{thm}[subsection]{Theorem}
\newtheorem{lem}[subsection]{Lemma}
\newtheorem{cor}[subsection]{Corollary}

\newcommand{\mbf}{\mathbf}
\newcommand{\mbb}{\mathbb}
\newcommand{\mrm}{\mathrm}
\newcommand{\bin}[2]{
\left[
   \begin{array}{@{}c@{}}
     #1 \\#2
   \end{array}
\right]    }     
\newcommand{\s}{\mbf s}
\newcommand{\mb}{\mbf b}
\newcommand{\mc}{\mbf c}
\newcommand{\md}{\mbf d}
\newcommand{\mg}{\mbf g}
\newcommand{\A}{\mathbb A}
\newcommand{\E}{\mathrm E}
\newcommand{\G}{\mathrm G}
\newcommand{\Q}{\mathbb Q}
\newcommand{\mD}{\mathcal D}
\newcommand{\mF}{\mathcal F}
\newcommand{\mP}{\mathcal P}
\newcommand{\mQ}{\mathcal Q}
\newcommand{\mS}{\mathcal{S}}
\newcommand{\tF}{\tilde{\mF}}
 
\newcommand{\Vb}{V^{\bullet}}
\newcommand{\Tb}{T^{\bullet}}
\newcommand{\Wb}{W^{\bullet}}

\newcommand{\Hom}{\mathrm{Hom}}
\newcommand{\Ind}{\mrm{Ind}}
\newcommand{\Aut}{\mathrm{Aut}}

\title[Notes on affine canonical and monomial bases]
{Notes on affine canonical and monomial bases} 
\author{ Yiqiang Li}
\address{Department of Mathematics\\ Yale University 
\\10 Hillhouse Avenue\\P.O. Box 208283\\ New Haven, CT 06520}
\email{yiqiang.li@yale.edu}
\thanks{}

\begin{document}

\begin{abstract}
We investigate the affine canonical basis 
(\cite{lusztig2}) and the monomial basis constructed in ~\cite{LXZ} in
Lusztig's geometric setting. We show that the transition matrix
between the two bases is upper triangular with 1's in the diagonal and
coefficients in the upper diagonal entries in 
$\mathbb Z_{\geq 0}[v, v^{-1}]$. As a consequence, we show that part
of the monomial basis elements give rise to resolutions of support
varieties of the affine canonical basis elements as simple perverse sheaves.
\end{abstract}

\maketitle

\section{Introduction}
\label{introduction}
Let $\mbf U^-$ be the negative part of the quantized enveloping
algebra 
associated to a  Cartan matrix $\mrm C$.
The canonical basis $\mbf B$ for $\mbf U^-$ was first constructed by
Lusztig in ~\cite{lusztig1} when $\mrm C$ is
symmetric and positive definite. 
In ~\cite{lusztig1}, Lusztig investigated the algebra $\mbf U^-$ in
three settings: algebraic, quiver and  geometric, 
where the canonical basis can be constructed from.
Here the terminology ``geometric'' means that 
the theory of perverse sheaves (\cite{BBD}) is used, 
``quiver'' means that we need to work on the setting of Ringel-Hall
algebras defined in ~\cite{Ringel}, 
and ``algebraic" means that it is purely Lie theoretic, 
i.e.,  no algebraic geometry or quiver representation results are used.
In fact, Lusztig first defined $\mbf B$ algebraically in
~\cite{lusztig1}. Then he studied $\mbf B$ in the framework of 
Ringel-Hall algebras and showed that $\mbf B$ is really a basis of $\mbf U^-$.
The geometric situation is finally brought in and produces various
remarkable properties  of $\mbf B$, for example, positivity and integrality.

The geometric approach to construct the canonical basis $\mbf B$ was
further studied in ~\cite{lusztig2} for $\mrm C$ symmetric  
and in ~\cite{lusztig5} for $\mrm C$ arbitrary.
The algebraic approach was further studied in ~\cite{BCP} for $\mrm C$ affine. 
(See ~\cite{Kashiwara} for a different (algebraic) approach, 
Kashiwara's approach works for arbitrary $\mrm C$.
Note that the canonical basis $\mbf B$ coincides with 
Kashiwara's global crystal basis in ~\cite{Kashiwara}, as was shown
in ~\cite{GL}.)
Recently, the quiver approach was carried out by Z. Lin, J. Xiao and
G. Zhang  in ~\cite{LXZ} for $\mrm C$ affine and symmetric. 
In ~\cite{LXZ}, a PBW-basis $E_{\mathcal M}$ is constructed 
by using  the representation theory of affine quivers.
The bar involution is upper triangular with the diagonal entries equal
to 1 with respect to the PBW-basis $E_{\mathcal M}$. By a method of
Lusztig (\cite{lusztig1}), one can deduce a unique bar invariant basis
$\mathcal{E_M}$ such that the transition matrix between 
$E_{\mathcal  M}$ and  $\mathcal{E_M}$ is upper triangular with
entries in the diagonal equal to $1$ and entries above the diagonal in
$v^{-1}\mathbb Q[v^{-1}]$.
In ~\cite{N}, Nakajima defined a modified PBW-basis $L_{\mathcal M}$
(of $E_{\mathcal M}$). The bar involution on  the basis
$L_{\mathcal M}$ has the same property as that on $E_{\mathcal M}$. 
Again by applying Lusztig's method on Nakajima's PBW-basis
$L_{\mathcal M}$, one still has a bar invariant basis. Nakajima
conjectured in ~\cite{N} that the latter bar invariant basis coincides with
Lusztig's affine canonical basis.

This paper is an attempt to understand the construction of 
the PBW-basis $E_{\mathcal M}$ and the bar invariant basis $\mathcal{E_M}$
in ~\cite{LXZ} in Lusztig's geometric setting. 
Unlike the situation of quivers of finite type, 
the geometric counterpart of elements in $E_{\mathcal M}$ are not in
Lusztig's algebra. 
Instead, we study the monomial basis $\mathcal F$ in ~\cite{LXZ} in the geometric
setting. It turns out that
the transition matrix $M$ between $\mathcal F$ and $\mbf B$ is
upper triangular, the entries in the diagonal are $1$ and
the entries above the diagonal are in $\mathbb Z_{\geq 0} [v,v^{-1}]$.
(This implies that the transition matrix between $\mbf B$ and
$\mathcal{E_M}$ is upper triangular with coefficients in 
$\mathbb Q[v,v^{-1}]$ and 1's in the diagonal entries.)
Moreover, certain monomials in $\mathcal F$ provide the
resolutions of the singularities of the support varieties of affine
canonical basis elements (regarded as simple perverse sheaves). 
  
In order to prove that the upper triangularity of the two bases
$\mathcal F$ and $\mbf B$, the representation-directed property of the
representation theory of affine quivers is essentially used. 
The proof is an adaption of Lusztig's argument in
~\cite{lusztig3}. It is somewhat simplified by using a lemma in ~\cite{Reineke}.

In Section ~\ref{affinequivers}, we review the theory of affine
quivers. 
We give a brief description of the classification of 
indecomposable representations of affine quivers.
In Section ~\ref{geometric}, we review Lusztig's geometric realization
of the canonical basis. 
In Section ~\ref{momomial}, we study some monomials in special orders,
which will be used as components in making the monomial basis $\mathcal F$. 
Then we identify the monomial basis $\mathcal F$ with the monomial
basis $\{\mbf{m_{c}}\;|\; \mbf c \in \mathcal M\}$ in ~\cite{LXZ}.
In the last section, we investigate the relationship between $\mathcal F$ and $\mbf B$. 
Finally, we draw the conclusions of the major results in this paper.

{\bf Acknowledgements.} We thank Professor Z. Lin 
for giving the preprints \cite{N} and ~\cite{LXZ} 
and suggesting the problem to the author.
The author is indebted to Professor H. Nakajima, who
pointed out several mistakes in this paper and explained in
details his conjecture in ~\cite{N} to the author. 
We also thank Professor J. Xiao for
informing the author his recent progress on Nakajima's conjecture.

\section{Representation theory of affine quivers}
\label{affinequivers}
We give a review of representation theory of affine quivers.

\subsection{Graphs, Root systems and Weyl groups}
\label{graph}
Let $I$ be a set. Denote by $\mbf I$ the set consisting of all two-element subsets of $I$. 
A $graph$ is a triple $\Gamma=(I,H, e: H\to \mbf I)$ where $I$ and $H$
are finite sets and $e$ is a map.
We call $I$  (resp. $H$) the $vertex$ (resp. $edge$) set of $\Gamma$. 

{\em
All graphs considered in this paper will be of affine type 
$A_n^{(1)}$ $(n\geq 2)$, $D_n^{(1)}$ $(n\geq 4)$ and $E_n^{(1)}$ $(n=6,7,8)$.
}

The $symmetric$ $Euler$ $form$ 
$(,): \mathbb Z[I]\times \mathbb Z[I]\to \mathbb Z$ associated to 
a given graph $\Gamma=(I,H,e)$ is defined by 
\begin{align}
(i, i)&=2 \quad \text{for all}\;  i\in I.\\
(i, j)&=- \#\{ h\in H\,|\, e(h)=\{i,j\}\} \quad \text{for all} \; i\neq j \in I.
\end{align}

A nonzero element $\alpha$ in $\mathbb Z[I]$ is called a $root$ if $(\alpha, \alpha) \leq 2$.
Denote by $R$ the set of all roots. Let $R_+=R\cap \mathbb Z_{>0}[I]$ and $R_-=-R_+$. 
It is well-known that 
\[
R=R_+ \sqcup R_-.
\]
One can check that $(\alpha,\alpha)\in 2\mathbb Z$. Thus if $\alpha$
is a root, then $(\alpha,\alpha)=0$ or $2$.
A root is called $real$ (resp. $imaginary$) if $(\alpha,\alpha)=2$ (resp. $0$). 
Denote by $R_+^{re}$ (resp. $R_+^{im}$) the set of all real (resp. imaginary) positive roots. 
It is well-known that all imaginary positive roots are of the form $n\delta$ where $n\in \mathbb Z_{>0}$ and $\delta \in R_+$.
$\delta$ is called the minimal imaginary positive root.
An {\em extending} vertex $i\in I$ is a vertex such that $\delta_i=1$.

Given any $i\in I$, we denote by $s_i$ the automorphism of $\mathbb Z[I]$ given by 
\[
s_i: \alpha \mapsto \alpha-(\alpha, i)i
\]
for any $\alpha \in \mathbb Z[I]$.
Note that $s_i$ is involutive, i.e., $s_i s_i=\text{id}_{\mathbb Z[I]}$. 
The bilinear form $(,)$ is $s_i$-invariant for all $i\in I$.
The $Weyl$ $group$ $W$ of $\Gamma$ is the subgroup of 
the automorphism group of $\mathbb Z[I]$ generated by $s_i$ for all $i\in I$.

\subsection{Quivers}
\label{quiver}

Let $\Gamma=(I,H,e)$ be a graph. An $orientation$ of $\Gamma$  is 
a pair of maps $s,t: H\to I$ such that $\{s(h),t(h)\}=e(h)$ for all $h\in H$.
In other words, the pair $(s,t)$ gives an orientation for each edge in $H$.

The graph $\Gamma$ equipped with the orientation $(s,t)$ will be called a $quiver$, denoted by $Q=(I,H,e,s,t)$.
We call $\Gamma$ the $underlying$ $graph$ of $Q$.
A quiver is called {\em affine} (or $tame$) if its underlying graph $\Gamma$ is affine.

{\em All quivers considered in this paper will be affine.}

Given any $h\in H$, we call $s(h)$ (resp. $t(h)$) the $starting$ (resp. $terminating$) vertex of the edge $h$. 
We also call edges in $H$ {\em arrows} in $Q$. 
Pictorially, we write $h: s(h) \to t(h)$ or $s(h) \overset{h}{\to} t(h)$ 
for indicating the starting and terminating vertex of $h$.

Note that once the pair $(s,t)$ is given, the map $e$ is then completely determined. 
We simply write $Q=(I,H, s, t)$ for a quiver.

An {\em oriented cycle} in $Q$ is a sequence of arrows $h_1, \cdots, h_m$ such that
$t(h_1)=s(h_2)$, $t(h_2)=s(h_3)$, $\cdots$, $t(h_{m-1})=s(h_m)$ and $t(h_m)=s(h_1)$.
A {\em cyclic quiver} is a quiver $Q=(I,H,s,t)$ of type $A_n^{(1)}$, where we order $I$ as $1, \cdots n$ 
and $H$ as $h_1, \cdots h_n$, such that $s(h_m)=m$ and $t(h_m)=m+1$ for all $m=1,\cdots n$ (we set $n+1=1$).
Since $Q$ is always affine in this paper, $Q$ has no oriented cycles except that $Q$ is a cyclic quiver.

Given any quiver $Q=(I,H,s,t)$, the $Euler$ $form$ of $Q$:
\[
<,>: \mathbb Z[I]\times \mathbb Z[I] \to \mathbb Z
\]
is defined by
\[
<\alpha,\beta>=\sum_{i\in I} \alpha_i \beta_i -\sum_{h\in H} \alpha_{s(h)}\beta_{t(h)}
\]
for any $\alpha=\sum_{i\in I}\alpha_i \ i,\beta=\sum_{i\in I} \beta_i\ i \in \mathbb Z[I]$. By definitions, we have
\begin{equation}
(\alpha,\beta)=<\alpha,\beta>+<\beta,\alpha> \quad \text{for any}\; \alpha, \beta \in \mathbb Z[I].
\end{equation}
Given any $i\neq j \in I$, we set
\[
a_{ij}= \# \{ h\in H \; | \; s(h)=i, t(h)=j \}.
\]

\subsection{Representation of affine quivers}
\label{representation}
Fix an algebraically closed field $k$.

A {\em representation} of $Q$ over $k$ is a pair $(V,x)$ where $V=\oplus_{i\in I} V_i$ is an $I$-graded $k$-vector space and 
$x$ is a collection of $k$-linear maps $x_h: V_{s(h)} \to V_{t(h)}$ for all $h\in H$.

For two given representations $\mbf V=(V,x)$ and $\mbf W=(W,y)$, 
a {\em morphism} $\phi$ between $\mbf V$ and $\mbf W$ is a collection
of linear maps $\phi_i : V_i\to W_i$ for all $i\in I$ 
such that $\phi_{t(h)} x_h=y_h \phi_{s(h)}$ for all $h\in H$.

These define an abelian category, denoted by $\mrm{Rep} (Q)$, whose objects are 
representations of $Q$ and morphisms are morphisms between 
representations.

A {\em nilpotent} representation $\mbf V$ is a representation satisfying the condition:
there exists $N$ such that for any $h_1,\cdots,h_N$ in $H$ 
satisfying $s(h_m)=t(h_{m-1})$ ($m=2,\cdots, N)$) then the composition 
\[
x_{h_N}\, x_{h_{N-1}}\, \cdots x_{h_1}: V_{s(h_1)} \to V_{t(h_N)}
\] 
is a zero map.
Denote by $\mrm{Nil}(Q)$ the full subcategory of $\mrm{Rep}(Q)$ consisting of all nilpotent representations. 
If $Q$ has no oriented cycles, $\mrm{Nil}(Q)$ and $\mrm{Rep}(Q)$ coincide. 

The $dimension$ of a representation $\mbf V=(V,x)$ of $Q$ is $\sum_{i\in I} \dim V_i\ i\in \mathbb N[I]$, 
denoted by $|\mbf V|$ or $|(V,x)|$. 
Given any  two representations $\mbf V$ and $\mbf W$ of $Q$, we write $\Hom_Q(\mbf V,\mbf W)$ 
and $\mrm{Ext}^1_Q(\mbf V,\mbf W)$
for 
$\Hom_{\mrm{Rep}(Q)}(\mbf V,\mbf W)$ and   $\mrm{Ext}^1_{\mrm{Rep}(Q)}(\mbf V,\mbf W)$,
respectively.

Given any representation $\mbf V=(V,x)$ and $\mbf W=(W,y)$, the following sequence is exact:
\[
\Hom_Q(\mbf V,\mbf W)\\
\overset{a}{\hookrightarrow} \oplus_{i\in I} \Hom_k(V_i,W_i) \overset{b}{\to} 
\oplus_{h\in H} \Hom(V_{s(h)}, W_{t(h)}) \overset{c}{\twoheadrightarrow} \mrm{Ext}^1_Q(\mbf V,\mbf W) 
\]
where 
$a: \phi \mapsto \oplus_{i\in I} \phi_i$, 
$b: \oplus_{i\in I} \phi_i \mapsto \oplus_{h\in H} (\phi_{t(h)} x_h-y_h \phi_{s(h)})$ and
\[
c: \oplus_{h\in H} \psi_h \mapsto (0\to \mbf V \to \mbf E \to \mbf W \to 0)
\]
with $\mbf E=\left(\oplus_{i\in I} V_i\oplus W_i, \oplus_{h\in H} \begin{pmatrix} x_h & 0 \\ \psi_h & y_h \end{pmatrix}\right)$. 
Here the hook (resp. double head) arrow $" \hookrightarrow "$ (resp. $"\twoheadrightarrow "$) represents that $a$ (resp. $b$) is injective
(resp. surjective). 

From the exact sequence above, we have
\begin{prop}
\label{here}
\begin{enumerate}
\item $<|\mbf V|, |\mbf W|>=\dim \Hom_Q(\mbf V,\mbf W)-\dim \mrm{Ext}^1_Q(\mbf V,\mbf W)$.
\item $\mrm{Rep}(Q)$ is hereditary, i.e., the extension groups $\mrm{Ext}^n_Q(\mbf V,\mbf W)$ vanish 
      for any $\mbf{ V, W}\in\mrm{Rep}(Q)$ and  $n\geq 2$.
\end{enumerate}
\end{prop}

{\bf Remark.} Proposition ~\ref{here} (1) justifies why the the form $<,>$ defined in Section ~\ref{quiver} is called Euler form.

\subsection{BGP-reflection functors}
\label{BGP}
A vertex $i \in I$ is called a $sink$ (resp., a $source$) if  
$i\in \{s(h),t(h)\}$ implies $t(h)=i$ (resp.,
$s(h)=i$) for any $h \in H$.

For any  $i \in I$,
let $\sigma_{i}Q=(I, H, s', t')$ be the quiver whose underlying graph is the same as $Q$ and  whose  orientation
$(s', t')$ is  defined by 
\begin{enumerate}
\item[] $s'(h)=s(h)$ and $t'(h)=t(h)$ for all $h\in H$ such that $i\notin \{s(h), t(h)\}$;

\item[] $s'(h)=t(h)$ and $t'(h)=s(h)$ for all $h\in H$ such that $i\in \{s(h),t(h)\}$.
\end{enumerate}
In other words, $\sigma_i Q$ is the quiver obtained by reversing 
the arrows in $Q$ that start or terminate at $i$.

Assume that $i$ is a sink. 
We set 
\[
H_i^+=\{ h\in H \; | \; t(h)=i\}.
\]
For any representation $(V,x)$ of $Q$, we set
\[
V_{H_i^+}=\oplus_{h\in H_i^+} V_{s(h)}
\quad \text{and} \quad 
x_i^+=\oplus_{h\in H_i^+} x_h: V_{H_i^+} \to V_i.
\]
We also denote by $\mrm{Ker} (x_i^+)$ (resp. $\mrm{Im} (x_i^+)$) the kernel (resp. image) of $x_i^+$. 
The {\em BGP-reflection functor} 
\[\Phi_{i}^+:\mrm{Rep}(Q) \to \mrm{Rep} (\sigma_{i}Q)\] 
with respect to $i$ is defined in the following way. 
For any $(V, x) \in \mrm{Rep}(Q)$,
$\Phi_{i}^+(V,x)=(W,y)\in \mrm{Rep}(\sigma_{i}Q)$, 
where $W=\oplus_{j\in I} \, W_j$ with
$W_j=V_j$ if $j\neq i$ and $W_i=\mrm{ker}(x_i^+)$ and 
$y=(y_h\; |\; h \in H)$ with
$y_h=x_h$ if $h\notin H_i^+$ and $y_h$ is the composition of the maps
$W_i \hookrightarrow V_{H_i^+} \twoheadrightarrow W_{t'(h)} \quad
\text{if} \; h\in H_i^+$.
Note that the assignments extend to a functor.

Similarly, assume that $i$ is a source. 
Let $H_i^-=\{h\in H \;|\; s(h)=i\}$. 
For any $(V,x)\in \mrm{Rep}(Q)$, we set
$V_{H_i^-}=\oplus_{h\in H_i^-} V_{t(h)}$
and
$x_i^-=\oplus_{h\in H_i^-} x_h: V_i \to V_{H_i^+}$.
Denote by $\mrm{Coker}(x_i^-)$ the cokernel of $x_i^-$.
The {\em BGP-reflection functor} 
\[\Phi_{i}^-:\mrm{Rep}(Q) \to \mrm{Rep} (\sigma_{i}Q)\] 
with respect to $i$ is defined in the following way. 
For any $(V, x) \in \mrm{Rep}(Q)$,
$\Phi_{i}^-(V,x)=(W,y)\in \mrm{Rep}(\sigma_{i}Q)$, 
where $W_j=V_j$ if $j\neq i$, $W_i=\mrm{Coker}(x_i^+)$, 
$y_h=x_h$ if $h\notin H_i^-$ and $y_h$ is the composition of the maps
$W_{s'(h)} \hookrightarrow V_{H_i^-} \twoheadrightarrow W_i \quad
\text{if} \; h\in H_i^+$.
Note that the assignments again extend to a functor.

Assume that $i$ is a sink. 
For a given representation $\mbf V=(V,x)$, 
let $\mbf V(i)=(W,y)$ be a representation of $Q$ such that 
$W_j=0$ for $j\neq i$ and $W_i=V_i/\mrm{Im}(x_i^+)$ and $y\equiv 0$.
From the construction of $\Phi_i^+$ and $\Phi_i^-$, 
one can deduce that (see ~\cite{BGP})
\begin{align}
\label{BGPlemma}
\mbf V \simeq \Phi_i^-\Phi_i^+ (\mbf V) \oplus \mbf V(i).
\end{align}
Denote by $\mrm{Rep}_i^+(Q)$ the full subcategory of $\mrm{Rep}(Q)$ whose objects are representations
such that $x_i^+$ is surjective. Similar, when $i$ is a source, denote by $\mrm{Rep}_i^-(Q)$ the full subcategory of
$\mrm{Rep}(Q)$ whose objects are representations such that $x_i^-$ is injective. 
From the definitions and (\ref{BGPlemma}), one has 
\begin{prop} Assume that $i$ is a sink, then
\begin{enumerate}
\item $\Phi_i^+(\mbf V) \in \mrm{Rep}_i^-(\sigma_iQ)$ for any $\mbf V\in \mrm{Rep}(Q)$.
\item $|\Phi_i^+(\mbf V)|=s_i(|\mbf V|)$ for any $\mbf V\in \mrm{Rep}_i^+(Q)$.
\item The restriction $\Phi_i^+: \mrm{Rep}_i^+(Q)\to \mrm{Rep}_i^-(\sigma_i Q)$ defines an equivalence of categories.
      Its inverse is $\Phi_i^-: \mrm{Rep}_i^-(\sigma_i Q)\to \mrm{Rep}_i^+(Q)$.
\end{enumerate}
\end{prop}

\subsection{Classification of indecomposable representations}
\label{classification}

Assume that $Q$ is not  a cyclic quiver.
Then $Q$ has no oriented cycles.
We can order the vertex set $I$ ($|I|=n+1$) of Q in a way, say  
\[
i_0, i_1 ,\cdots,i_n
\] 
such  that
$i_0$ is a sink of $Q$, $i_1$ is a sink of $\sigma_{i_0}Q$, $\cdots$,
$i_j$ is a sink of $\sigma_{i_{j-1}}\cdots \sigma_{i_0}Q$, $\cdots$.
In other words, 
$i_n$ is a source of $Q$, $i_{n-1}$ is a source of $\sigma_{i_n}Q$, $\cdots$,
$i_j$ is a source of $\sigma_{i_{j+1}} \cdots \sigma_{i_n}Q$, $\cdots$.
Since there is no oriented cycles in $Q$, one can show by induction (see ~\cite{BGP}) that 
\[
Q=\sigma_{i_n}\cdots \sigma_{i_0} Q
\quad 
\text{and}
\quad
Q=\sigma_{i_0}\cdots \sigma_{i_n} Q.
\] 
Define the Coxeter functor $\Phi^+: \mrm{Rep} (Q) \to
\mrm{Rep} (Q)$ (resp. $\Phi^-: \mrm{Rep}(Q) \to\mrm{Rep}(Q)$)
to be the composition of functors
\[
\Phi^+ = \Phi_{i_n}^+\circ \cdots \circ \Phi_{i_1}^+\Phi_{i_0}^+\quad
\text{(resp.}\; 
\Phi^-=\Phi_{i_0}^-\Phi_{i_1}^- \circ \cdots \circ \Phi_{i_{n}}^-).
\]

Denote by $\Ind (Q)$ the set of all isomorphism classes of 
indecomposable representations in $\mrm{Rep}(Q)$.
For any representation $\mbf V$, we write $[\mbf V]$ for its isomorphism class. 
We identify $\Ind(Q)$ with the set of representatives of 
indecomposable representations in $\mrm{Rep}(Q)$.
By abusing of notation, we write $\mbf V \in \Ind (Q)$ for $[\mbf V]\in \Ind (Q)$.

An indecomposable representation $\mbf V$ is called {\em preprojective} if
$(\Phi^+)^m(\mbf V)=0$ for $m$ large enough;
{\em preinjective} if $(\Phi^-)^m (\mbf V)=0$ for $m$ large enough and
{\em regular} if $(\Phi^+)^m(\mbf V)\neq 0$ for arbitrary $m$.

A regular indecomposable representation $\mbf V$ is called
{\em homogeneous} if $\Phi^+(\mbf V)\simeq \mbf V$ and
{\em inhomogeneous} if $\Phi^+(\mbf V)\not \simeq \mbf V$.

A representation $\mbf V$ is called preprojective 
if all the indecomposable summands of $\mbf V$ are preprojective. 
We define  a representation to be  
preinjective, regular, homogeneous regular and inhomogeneous regular in a similar way.

For any inhomogeneous regular indecomposable representation $\mbf V$, there
exists $p$ such that $(\Phi^+)^p(\mbf V)\simeq \mbf V$ (see ~\cite{BGP}).
The {\em period} of an inhomogeneous regular representation $\mbf V$
is the smallest positive integer
$p$ such that $(\Phi^+)^p(\mbf V)\simeq \mbf V$.

Let $S_i$ be the representation corresponding to the vertex $i$. 
It's a representation $(V, x)$ where $V_i=k$, $V_j=0$  if $j\neq i$, and $x\equiv 0$. 
Clearly, $S_i$ is a simple object in $\mrm{Rep}(Q)$. 
Moreover, the set $\{S_i\;|\; i\in I\}$ is 
a complete list of pairwise non isomorphic simple representations in $\mrm{Rep}(Q)$.
Note that given a graph, the definition of the simple representation 
works for any orientation of the graph. 
By abuse of notation, we always denote by $S_{i}$ 
the simple representation corresponding to the vertex $i$ regardless of the orientation. 

We set
\[
P_0=S_{i_0}, P_1=\Phi_{i_0}^-(S_{i_1}), \cdots , P_n
=\Phi_{i_0}^-\Phi_{i_1}^-\cdots\Phi_{i_{n-1}}^-(S_{i_n}).
\]
Then the set  $\{P_m\;|\; m=0,\cdots,n\}$ is 
a complete set of pairwise non isomorphic projective representations of $\mrm{Rep}(Q)$.
For any $m>n$, $m$ can be written as $m=(n+1)p+q$ 
with $p$ and $q$ positive integers and $0\leq q\leq n$.
We set
\[
P_m=(\Phi^-)^p (P_q).
\]
Then the set 
\[\mathscr P=\{P_m\;|\: m\in \mbb Z_{\geq 0}\}\] 
forms a complete set of 
pairwise non isomorphic preprojective representations of $\mrm{Rep}(Q)$.

Similarly, we set 
\[
I_n=S_{i_n}, I_{n-1}=\Phi_{i_n}^+(S_{i_{n-1}}), \cdots, I_0
=\Phi_{i_n}^+\cdots \Phi_{i_1}^+\Phi_{i_0}^+ (S_{i_0}).
\]
Then the set $\{I_m\;|\; m=0,\cdots,n\}$ is a complete set of pairwise non isomorphic injective representations of $\mrm{Rep}(Q)$.
For any $m<0$, $m$ can be written as $m=-((n+1)p+q)$ with $p\in \mathbb Z_{\geq 0}$ and $q$ positive integers and $0\leq q\leq n$.
We set
\[
I_m=(\Phi^+)^p (I_q).
\]
Then the set 
\[\mathscr I=\{I_m\;|\: m\in \mbb Z_{\leq n}\}\] 
forms a complete set of pairwise non isomorphic preinjective representations of $\mrm{Rep}(Q)$.

Let $\mrm{Reg}(Q)$ be the full subcategory of $\mrm{Rep}(Q)$
whose objects are regular representations. Then 
$\mrm{Reg}(Q)$ is an extension-closed  full subcategory of  $\mrm{Rep}(Q)$. Moreover,
the restriction of $\Phi^+$ on $\mrm{Reg}(Q)$ is an equivalence of
categories. $\Phi^-$ is its inverse.

The simple objects in $\mrm{Reg}(Q)$ are called {\em regular} {\em simple} representations. 
For each regular simple representation $T$ and $m\geq 1$, there is a
unique (up to an isomorphism) 
regular indecomposable representation, denoted by $R_{T,m}$, 
such that the composition factors of $R_{T,m}$
are $T$, $\Phi^+(T)$, $\cdots, (\Phi^+)^m (T)$.
Moreover, all regular indecomposable representations are obtained this way (see ~\cite{DR}).

Assume that $T$ is a  regular simple representation of period $p$.
The set 
\[
\{T,\cdots,(\Phi^+)^{p-1} \, T\}
\] 
is called the $\Phi^+$-{\em orbit} of $T$.
Furthermore, 
\[
|T|+|\Phi^+(T)|+\cdots+|\Phi^+)^{p-1}(T)|=\delta.
\]
Note that $|(\Phi^+)^r(T)|$ is a positive root of finite type, thus
$(\Phi^+)^r (T)$ has no self-extension for $r=0,\cdots, p-1$ (\cite{FMV}).
The {\em tube} $\mathscr T=\mathscr T_T$
is the set of all indecomposable regular representations whose 
regular composition factors belong to this orbit. 
Any regular indecomposable representation  belongs to a tube. 
Any representation in a tube has the same period, 
which is called the period of the tube.
All but finitely many tubes have period one.  
Tubes with period one will be called {\em homogeneous} tubes and tubes with
period $>1$ will be called {\em inhomogeneous} tubes. 
In fact, the number of inhomogeneous tubes is $\leq 3$ (\cite{DR}).
We denote them by $\mathscr T_1,\cdots, \mathscr T_s$ ($s\leq 3$). 
The indecomposable representations in $\mathscr T_i$  will be denoted by
$T_{i, a, l}$ where $a=1, \cdots p_i$ and $l\in \mathbb Z_{\geq 0}$. 
Here $p_i$ is the period of $\mathscr T_i$.
Their relations  are governed by $T_{i, a+1,l}=\Phi^+(T_{i,a,l})$ 
for $a=1,\cdots, p_i$ ($p_i+1:=1$) and $l\in \mathbb Z_{\geq 0}$.

Let $\mrm{Rep}(\mathscr T)$ be the full subcategory of $\mrm{Rep} (Q)$ 
whose objects are direct sums of indecomposable representations in
$\mathscr T$.

Given any representation $M\in \mrm{Rep}(\mathscr T_i)$ with $\mathscr T_i$ of period $p_i$, 
$M$ is called {\em aperiodic} if for any 
$l\in \mathbb N$,  not all the  representations 
\[
T_{i, 1, l}, T_{i,2, l}, \cdots, T_{i, p_i, l}
\]
are direct summands of $M$. 
From the above analysis, one can show 
\begin{lem} 
\label{vanishing}
Let $\mbf V$ and $ \mbf W\in \mrm{Ind}(Q)$ be one of the following cases.
\begin{itemize}
\item $\mbf V=I_m$ and $ \mbf W=I_{m'}$ for $m\leq m'\in\mathbb Z_{\leq n}$.
\item $\mbf V=P_m$ and $\mbf W=P_{m'}$ for $m\leq m'\in\mathbb Z_{\geq 0}$.
\item $\mbf V$ and $\mbf W$ are both regular, but they are in different tubes.
\item $\mbf V$ is non preinjective and $\mbf W$ is preinjective.
\item $\mbf V$ is preprojective and $\mbf W$ is non preprojective. 
\end{itemize}
Then we have
\begin{align}
&\mrm{Ext}^1_Q(\mbf V,\mbf W)=0.\\
&\mrm{Hom}_Q(\mbf V, \mbf W)=0 \quad \text{if} \; \mbf V \; \text{and}
\; \mbf W\;\text{are not isomorphic}.
\end{align}
\end{lem}
 
The following lemma is proved in ~\cite{Reineke}.
\begin{lem}
\label{Reineke}
Assume that $\mbf V(1)=(V(1),x(1)), \cdots, \mbf V(s)=(V(s),x(s))$ in
$\mrm{Rep}(Q)$ satisfy $\mrm{Ext}^1_Q(\mbf V(k), \mbf V(l))=0$ for all
$k < l$. Suppose that 
\[
\mbf V^{\bullet}= (V=V^0 \supseteq V^1
\supseteq \cdots \supseteq V^{s}=0)
\]
is a flag such that
$|V^{l-1}/V^{l}|=|V(l)|$ for all $l=1,\cdots, s$. If $x\in \E_V$
stabilizes $\mbf V^{\bullet}$ and the restriction $x|_{V^{l-1}/V^l}$
of $x$ to $V^{l-1}/V^l$ is in the closure of the $\G_{V(l)}$-orbit
$O_{x(l)}$ of $x(l)$ for $l=1, \cdots s$, then $x$ is in the closure
of $x(1)\oplus \cdots \oplus x(s)$. 

Moreover, if $\mrm{Hom}_Q(V(k), V(l))=0$ for $k >l$ and 
$(V,x) \simeq (V, x(1)\oplus \cdots \oplus x(s))$, then
$x|_{V^{l-1}/V^l}$ is in $O_{x(l)}$ for any $l=1,\cdots, s$ and $x$
fixes a unique flag of type $(|V(1)|, \cdots, |V(s)|)$. 
\end{lem}

\section{Lusztig's geometric realization of the canonical basis}
\label{geometric}

We give a review of Lusztig's geometric construction of the canonical basis.

\subsection{Representation spaces of $Q$}

Let $k$ be a fixed algebraically closed field. 
Let $V=\oplus_{i\in I} \, V_i$ be an $I$-graded $k$ vector space.
Let $|V|$ be its dimension $\sum_{i\in I} \, \dim V_i \, i\in \mbb N[I]$.
We set
\begin{align}
\E_V=\oplus_{h\in H} \, \Hom(V_{s(h)}, V_{t(h)}) 
\quad \text{and} \quad 
\G_V=\oplus_{i\in I} \Aut(V_i)
\end{align}
where $\Aut(V_i)$ is the group of all linear isomorphism of $V_i$.

For any $g=(g_i\; |\; i\in I) \in \G_V$, $x=(x_h\; |\; h\in H)\in
\E_V$, define $g.x=((g.x)_h\; |\; h\in h)$ by
\[
(g.x)_h= g_{t(h)}\,   x_h \,  g_{s(h)}^{-1}
\]
for any $h\in H$.
This defines a $\G_V$ action on $\E_V$.

Note that given any $x\in \E_V$, the pair $(V,x)$ is a representation of $Q$. The isomorphism classes of representations of 
$Q$ of dimension $|V|$ is then in one-to-one correspondence with the $\G_V$-orbits in $\E_V$.

\subsection{Flag varieties}
\label{flag}
Given any $\nu \in \mbb N[I]$, we define the set $\mS_{\nu}$ to be the
set consists of all sequences
\[
\s=(s_1 \, i_1, \cdots, s_n\, i_n)
\]
such that $\sum_{m} s_m \, i_m =\nu$.

Fix a $V$ such that $|V|=\nu$. For any $\s \in \mS_{\nu}$, we say that
a flag
\[
\Vb=(V=V^0\supseteq V^1 \supseteq \cdots \supseteq V^n=0)
\]
is of type $\s$ if $|V^{m-1}/V^m|=s_m\, i_m$, for all $m=1,\cdots n$.

we define the variety $\mF_{\s}$ to be
the variety consisting of all flags of type $\s$. 
Note that $\G_V$ acts transitively on $\mF_{\s}$.

Given $x \in \E_V$, $\Vb \in \mF_{\s}$, we say that $V^{\bullet}$ is $x$-stable if 
$x_h(V^m_{s(h)}) \subseteq V^m_{t(h)}$ for any $h\in H$ and
$m=1,\cdots, n$.

We define $\tF_{\s}$ to be the variety consisting of all pairs 
$(x, \Vb) \in \E_V \times \mF_{\s}$ such that $\Vb$ is
$x$-stable. Note that $\G_V$ acts on $\tF_{\s}$ naturally.
We have

\begin{lem}
\label{flagvarieties}
\begin{enumerate}

\item The variety $\mF_{\s}$ is smooth, irreducible, projective
      variety of dimension
      \[\dim \mF_{\s} = \sum_{m>m': \, i_m=i_{m'}}\, s_m \, s_{m'}.\]

\item The variety $\tF_{\s}$ is smooth, irreducible variety of
      dimension
     \[\dim \tF_{\s}=\sum_{m>m': \, i_m=i_{m'}}\, s_m \, s_{m'} +
       \sum_{h\in H:\, m\leq m'} \, a_{ij} \, s_m \, s_{m'}.\]

\item The first projection $\pi_{\s}: \tF_{\s} \to \E_V \quad ((x, \Vb)
      \mapsto x)$ is a $\G_V$-equivariant, proper morphism.
\end{enumerate}
\end{lem}

For a proof, see ~\cite{lusztig2}.

\subsection{Notations}
We fix some notations, most of them are taken from  ~\cite{lusztig5}.

Fix a prime $l$ that is invertible in $k$.
Given any algebraic variety $X$ over $k$, 
denote by $\mathcal{D}(X)$ the bounded derived category of complexes 
of $l$-adic sheaves on $X$ (\cite{BBD}).
Let $\mathcal M(X)$ be the full subcategory of $\mathcal D(X)$
consisting of all perverse sheaves on $X$ (\cite{BBD}).

Let $G$ be a connected algebraic group. Assume that $G$ acts on $X$
algebraically. Denote by $\mathcal D_G(X)$ the full subcategory of
$\mathcal D(X)$ consisting of all $G$-equivariant complexes over
$X$.
Similarly, denote by $\mathcal M_G(X)$ the full subcategory of
$\mathcal M(X)$ consisting of all $G$-equivariant perverse sheaves (\cite{lusztig5}).
Let $Y$ be a smooth, locally closed, irreducible $G$-invariant subvariety of $X$ and 
$\mathcal L$ an irreducible, $G$-equivariant, local system on $Y$.
Denote by $j: Y\to X$ the natural embedding. 
From \cite{BBD} and  ~\cite{BL},
the complex 
\[\mrm{IC}(Y, \mathcal L):=j_{!\star}(\mathcal L)[\dim Y]\] 
is a simple $G$-equivariant perverse sheaf on $X$.
Moreover, all simple $G$-equivariant perverse sheaves on $X$ are of this form.

Let $\bar{\mathbb Q}_l$ be an algebraic closure of the field of $l$-adic numbers.
By abuse of notation, denote by  $\bar{\mathbb Q}_l=(\bar{\mathbb Q}_l)_X$ 
the complex concentrated on degree zero, corresponding to
the constant $l$-adic sheaf over $X$. 
For any complex $K \in \mathcal D(X)$ and $n\in \mathbb Z$, 
let $K[n]$ be the complex such that $K[n]^i=K^{n+i}$ and the
differential is multiplied by a factor $(-1)^n$.
Denote by $\mathcal M(X)[n]$ the full subcategory of 
$\mathcal D(X)$ whose objects are of the form $K[n]$ with 
$K\in \mathcal M(X)$.
For any $K\in \mathcal D(X)$ and $L\in \mathcal D(Y)$, denote
by $K\boxtimes L$ the external tensor product of $K$ and $L$ in
$\mathcal D(X\times Y)$.

Let $f: X\to Y$ be a morphism of varieties, denote by
$f^*: \mathcal D(Y) \to \mathcal D(X)$ and 
$f_!: \mathcal D(X) \to \mathcal D(Y)$ 
the inverse image functor and the direct image functor with compact support, respectively.

If $G$ acts on $X$ algebraically and $f$ is a principal $G$-bundle,
then $f^*$ induces a functor (still denote by $f^*$) of equivalence
between $\mathcal M(Y)[\dim G]$ and $\mathcal M_G(X)$.
Its inverse functor is denoted by $f_{\flat}: \mathcal M_G(X) \to
\mathcal M(Y)[\dim G]$ (\cite{lusztig5}).

\subsection{Semisimple complexes on $\E_V$}
\label{semisimple}
Let $\bar{\Q}_l$ be the constant sheaf on $\tF_{\s}$. 
By Lemma ~\ref{flagvarieties} (2), the complex $\bar{\Q}_l\, [\dim
\tF_{\s}]$ is a $\G_V$-equivariant simple perverse sheaf on $\tF_{\s}$.
So by Lemma ~\ref{flagvarieties} (3), the complex
\[
L_{\s}:= (\pi_{\s})_!\, (\bar{\Q}_l \, [\dim \tF_{\s}])
\]
is  a $\G_V$-equivariant semisimple complex on $\E_V$.
\begin{example}
\label{example}
Let $V$ be an $I$-graded vector space over $k$ of dimension $mi\in \mathbb N[I]$. 
By definition, $\E_V$ contains a single element $0$. Let $\s=(ni)$. Then
$L_{\s}$ is the constant sheaf  on $\E_V$. We write $F(ni)$ for $L_{\s}$. 
\end{example}

Let $\mP_V$ be the set of all isomorphism classes of simple perverse
sheaves on $\E_V$ appearing as direct summands with possible shifts in $L_{\s}$, for all
$\s \in \mS_{\nu}$. Here $|V|=\nu$.

Let $\mQ_V$ be the full subcategory of $\mD(\E_V)$ whose objects
are finite direct sum of shifts of  simple perverse sheaves coming from $\mP_V$.
Note that all complexes in $\mQ_V$ are semisimple and $\G_V$-equivariant.

Let $\mQ(\E_T\times \E_W)$ be the full subcategory of $\mD(\E_T\times
\E_W)$ whose objects are finite direct sum of shifts of the simple
perverse sheaves of the form $K\boxtimes L$ for all $K\in\mQ_T$ and
$L\in \mQ_W$.

\subsection{Lusztig's induction functors}
\label{induction}
Let $W \subseteq V$ be an $I$-graded subspace. 
Let $T=V/W$ and $p:V\to T$ the natural projection. We sometimes write
$p_T$ to avoid confusion.

Given any $x\in\E_V$,
$W$ is called $x$-stable  if 
$x_h(W_{s(h)})\subseteq W_{t(h)}$, for all $h\in H$.

If $W$ is $x$-stable, it induces two elements $x_W$ and $x_T$ in $\E_W$ and
$\E_T$ respectively as follows. $(x_W)_h$ is the restriction of
$x_h$ to $W$, for all $h\in H$. $(x_T)_h$ is defined such that
$p_{t(h)}\, x_h=(x_T)_h\, p_{s(h)}$ for all $h\in H$.

We consider the following diagram
\begin{equation*}
\label{ind}
\E_T\times \E_W \overset{q_1}{\longleftarrow} \E'
\overset{q_2}{\longrightarrow}
 \E'' \overset{q_3}{\longrightarrow} \E_V,
\tag{*}
\end{equation*}
where
$\E''= \{ (x,V')\; |\; V' \; \text{is} \; x \text{-stable}, |V'|=|W|\}$;

$\E'$ is the variety consisting of all quadruples $ (x, V', r', r'')$ such that
$(x,V')$ is in  $\E''$,  $r': V/V' \to T$ and  $r'':V' \to W$ are graded linear isomorphisms.

$q_3: (x,V') \mapsto x$, $q_2: (x, V', r', r'') \mapsto (x, V')$
and $q_1: (x, V', r', r'') \mapsto (y', y'')$ with 
$y'_h=r'_{t(h)}\,  (x_{V/V'})_h \, (r')^{-1}_{s(h)}$ and 
$y''_h=r''_{t(h)}\, (x_{V'})_h \, (r'')^{-1}_{s(h)}$, for all $h\in H$.

(For convenience, we write $y'= r' \cdot (x_{V/V'}) \cdot (r')^{-1}$
and $y''= r'' \cdot (x_{V'}) \cdot (r'')^{-1}$.)

By definition, we have 
\begin{enumerate}

\item[(a)] $q_3$ is proper;

\item[(b)] $q_2$ is $\G_{T}\times \G_{W}$-principal bundle of fiber
           dimension 

           $d_2=\sum_{i\in I} (|T_i|^2+ |W_i|^2)$;
\item[(c)] $q_1$ is smooth with connected fibers of fiber dimension

      $d_1=\sum_{i\in I} (|T_i|^2+|W_i|^2) +\sum_{h\in H}
      |T_{s(h)}|\,|W_{t(h)}| +\sum_{i\in I} |T_i|\,|W_i|$.
\end{enumerate}

From diagram (\ref{ind}) and the above properties, 
we have a functor 
\[
(q_{3})_!\; (q_{2})_{\flat}\; q_1^*:\mQ(\E_T\times \E_W) \to \mD(\E_V).
\]
Given any $K\in \mQ_T, L\in \mQ_W$, we set
\[
K\circ L= (q_{3})_!\; (q_{2})_{\flat}\; q_1^* (K \boxtimes L) \, [d_1-d_2].
\]
We have
\begin{lem}
\label{ind1}
\begin{enumerate}

\item $K \circ L \in \mQ_V$.

\item $L_{\s'} \circ L_{\s''}= L_{\s'\s''}$, 
      where $\s'\in \mathcal S_{\tau}$ [$\tau=|T|$], 
      $\s'' \in \mathcal S_{\omega}$ ($\omega=|W|$), $\s'\s'' \in
      \mathcal S_{\nu}$ ($\nu=\tau + \omega$) is the
      concatenation of the two sequences.
\end{enumerate}
\end{lem}

\begin{proof}
We show that (2) holds first.
Consider the following diagram
\[
\begin{CD}
\tF_{\s'}\times \tF_{\s''} @<\tilde{q_1}<< \tilde{\E}' @>\tilde{q}_2>>
\tF_{\s'\s''}\\
@V\pi_{\s'}\times\pi_{\s''}VV  @V\pi'VV @V\pi''VV\\
\E_{T}\times \E_{W} @<q_1<< \E' @>q_2>> E'' @>q_3>> E_V,
\end{CD}
\]
where 
\begin{enumerate}
\item[(a)] $\tilde{\E}'$ is the variety consisting of all quadruples
      $(x, \Vb, r', r'')$ such that $(x, \Vb)$ is  in $ \tF_{\s'\s''}$ and
      $r': V/V^{l_0} \to T$, $r'': V^{l_0}\to W$ are linear
      isomorphisms, where $V^{l_0} \in \Vb$ satisfying
      $|V^{l_0}|=|W|$, $|V/V^{l_0}|=|T|$.

\item[(b)] $\tilde{q}_1: (x, \Vb, r', r'') \mapsto (y', \Tb), (y'', \Wb)$,
      where $(y', y'')$ is defined as in the definition of $q_1$,
      $\Wb=r''(\Vb \cap V^{l_0})$ and $\Tb=r'(p_{V/V^{l_0}}(\Vb))$.

\item[(c)] $\tilde{q}_2: (x, \Vb, r', r'') \mapsto (x, \Vb)$,  $\pi': (x, \Vb,
      r', r'')\mapsto (x, V^{l_0}, r', r'')$ and $\pi'': (x, \Vb) \mapsto
      (x, V^{l_0})$.
\end{enumerate}
One can check that the above squares are Cartesian. So we have
\begin{enumerate}
\item[(d)] $q_1^*\; (\pi_{\s'}\times \pi_{\s''})_*=\pi'_{*}
      \; (\tilde{q}_1)^*$.

\item[(e)] $\pi''_{*} \; (\tilde{q}_2)^*=(q_2)^* \; \pi''_{*}$.

\item[(f)] $\dim \tF_{\s'} \times \tF_{\s''} +d_1=\dim \tilde{\E}'=
      \dim \tF_{\s'\s''}+d_2$.
\end{enumerate}
Then
\begin{equation*}
\begin{split}
L_{\s'}\circ L_{\s''}
&=(q_3)_{!}\; (q_2)_{\flat}\; q_1^*\;
(L_{\s'}\boxtimes L_{\s''})[d_1-d_2]\\
&\overset{(\text{d})}{=}(q_{3})_!\; (q_2)_{\flat}\; \pi'_!(\bar{\Q}_l) [\dim
\tF_{\s'} +\dim \tF_{\s''}+d_1-d_2]\\
&\overset{(\text{e})}{=} (q_3)_!\; (q_2)_{\flat} \; q_2^* \; \pi''_!(\bar{\Q}_l)[\dim
\tF_{\s'} +\dim \tF_{\s''}+d_1-d_2]\\
&=(q_3)_!
\pi''_!(\bar{\Q}_l)[\dim \tF_{\s'} +\dim \tF_{\s''}+d_1-d_2]\\
&=(\pi_{\s'\s''})_! (\bar{\Q}_l)[\dim
\tF_{\s'} +\dim \tF_{\s''}+d_1-d_2]\\
&\overset{(\text{f})}{=}(\pi_{\s'\s''})_! (\bar{\Q}_l)[\dim \tF_{\s'\s''}]=L_{\s'\s''}
\end{split}
\end{equation*}
So (2) follows.

(1) follows from (2).
\end{proof}

Given any $K \in \mQ_X, L\in \mQ_Y$ and $ M \in \mQ_Z$, we have 

\begin{lem} [Associativity]
\label{associativity}
$(K\circ L) \circ M =K \circ (L\circ M)$.
\end{lem}

\begin{proof}
Assume that $X\oplus Y= U$, $U \oplus Z=V$.
Consider the commutative diagram
\[
\begin{CD}
(\E_X \times \E_Y) \times \E_Z @<t<< D\\
@Ap_1AA @Ao_1AA\\
\E'_{XY}\times \E_Z @<s<< C @>o_1>> D\\
@Vp_2VV @Vo_2VV @V n_2VV\\
\E''_{XY}\times \E_Z @<r_1<< B  @>r_2>> A\\
@Vp_3VV @Vo_3VV @Vn_3 VV\\
\E_U \times \E_Z @<q_1<< \E' @>q_2>> \E'' @>q_3>> \E_V,
\end{CD}
\]
where 
the bottom row is the diagram (\ref{ind}) with $T, W$ replaced
by $U, Z$, respectively.
The column to the left is the diagram (\ref{ind}) times $\E_Z$ with
$T, W$ replaced by $X, Y$. To avoid confusion, we write 
$\E'_{XY}, \E''_{XY}$ for $\E', \E''$.

$A$ is the variety consisting of all pairs $(x, \Vb)$,
where $x\in \E_V$, $\Vb=(V^1 \supseteq V^2)$ is
a flag in $V$, 
such that $|V/V^1|=|X|, |V^1/V^2|=|Y|, |V^2|=|Z|$ and $\Vb$ is $x$-stable.

$B$ is the variety consisting of all triples $(x, \Vb, \mb)$, 
where $(x, \Vb) \in A$, 
and $\mb=(b': V/V^2 \to U, b'': V^2 \to Z)$  is a
pair of linear isomorphisms.

$C$ is the variety consisting of all triples $(x, \Vb, \mb, \mc)$ where
$(x,\Vb, \mb) \in B$, 
$\mc=( c': V/V^1 \to X, c'': V^1/V^2 \to Y)$ is a
pair of linear isomorphisms.

$D$ is the variety consisting of all triples $(x, \Vb, \md)$, where
$(x, \Vb)$ are in $B$, $\md=(d': V/V^1 \to X, d'': V^1/V^2 \to Y, d''':V^2
\to Z)$ is a triple of linear isomorphisms.

The maps $q_1, q_2, q_3$ are from diagram (\ref{ind}) 
and $p_1, p_2, p_3$ are induced from diagram
(\ref{ind}) (by multiplying with the identity map $\text{id}: \E_Z \to \E_Z$).

$n_3: (x, \Vb) \mapsto (X, V^2)$, 
$o_3: (x, \Vb, \mb) \mapsto (x, V^2, \mb)$.

$r_2: (x, \Vb, \mb) \mapsto (x,\Vb)$, 

$r_1: (x, \Vb, \mb) \mapsto (( b' \cdot (x_{V/V^2}) \cdot (b')^{-1},
b'(V^1/V^2));  b'' \cdot (x_{V^2}) \cdot  (b'')^{-1})$.

$n_2: (x, \Vb, \md) \mapsto (x, \Vb)$,
$o_2: (x, \Vb, \mb, \md) \mapsto (x, \Vb, \mb)$,

$o_1: (x, \Vb, \mb, \mc) \mapsto (x, \Vb, \md)$,
where $\md=(c', c'', b'')$.

$s: (x, \Vb, \mb, \mc) \mapsto 
( (b' \cdot (x_{V/V^2}) \cdot (b')^{-1},
b'(V^1/V^2)), r', r'');  b'' \cdot (x_{V^2}) \cdot  (b'')^{-1}))$,
where $r'= c''\, (b')^{-1}, r''=c'\, (b')^{-1}$.

$t:(x, \Vb, \md) \mapsto 
(d'\cdot x_{V/V^1} \cdot (d')^{-1}, d'' \cdot x_{V^1/V^2} \cdot
(d'')^{-1}, d''' \cdot x_{V^2} \cdot (d''')^{-1}$.

$o_1$ is a $\G_U$-principal bundle, $o_2$ is a $\G_X\times
\G_Y$-principal bundle, $r_2$ is a $\G_U \times \G_Z$-principal
bundle, $n_2$ is a $\G_X\times \G_Y \times \G_Z$-principal bundle.

Since $r_2 \, o_2= n_2\, o_1$ and all maps are principal bundles, we
have $(r_2)_{\flat} \, (o_2)_{\flat} =(n_2)_{\flat}\, (o_1)_{\flat}$. 
So $(r_2)_{\flat} \, (o_2)_{\flat}\, o_1^*=(n_2)_{\flat}\,
(o_1)_{\flat}\, o_1^*= (n_2)_{\flat}$. I.e.
\begin{equation}
\label{n_2}
(r_2)_{\flat} \, (o_2)_{\flat}\, o_1^*=(n_2)_{\flat} 
\end{equation}

We denote by $d_1, d_2, e_1, e_2$ the fiber dimensions of the maps
$q_1, q_2, p_1, p_2$, respectively. We set $d=d_1-d_2+e_1-e_2$.

From the definitions, we have
\[
(K\circ L)\boxtimes M= (p_3)_! \, (p_2)_{\flat}\,
p_1^* (K\boxtimes L \boxtimes M)[e_1-e_2].
\]
For simplicity, we write $\mathcal M=K\boxtimes L \boxtimes M$.
So we have

\begin{equation*}
\begin{split}
(K\circ L) \circ M 
&=(q_3)_!\, (q_2)_{\flat}\, (q_1)^* \, (p_3)_! \, (p_2)_{\flat}\,
p_1^* (\mathcal M)\, [d]\\
&=(q_3)_!\, (q_2)_{\flat}\, (o_3)_!\, (r_1)^*\, (p_2)_{\flat}\, 
p_1^* (\mathcal M)\, [d]\\
&=(q_3)_!\, (n_3)_!\, (r_2)_{\flat}\, r_1^*\, (p_2)_{\flat}\, 
p_1^* (\mathcal M)\, [d]\\
&=(q_3)_!\, (n_3)_!\, (r_2)_{\flat}\, (o_2)_{\flat}\, s^* \,
p_1^* (\mathcal M)\, [d]\\
&= (q_3\, n_3)_! \, (r_2)_{\flat}\, (o_2)_{\flat}\,o_1^* \, t^*
(\mathcal M)\, [d]\\
&\overset{\ref{n_2}}{=} (q_3\, n_3)_! \,(n_2)_{\flat} t^*
(\mathcal M)\, [d].
\end{split}
\end{equation*}
The second to the forth equation  holds, since the squares
\[
(o_3, q_1, r_1, p_3), 
(o_3, q_2, r_2, n_3), 
\;\text{and}\; (o_2, r_1, s, p_2)\]
are Cartesian. 

Similarly, we can show that 
\[
K\circ (L\circ M)=
(q_3\, n_3)_! \,(n_2)_{\flat} t^* (\mathcal M)\, [d].
\]
Lemma follows.
\end{proof}

By Lemma ~\ref{associativity}, it will not cause any confusion when we
write $K \circ L \circ M$.

From the Proof above, we obtain a diagram
\[
\label{**}
\begin{CD}
\E_X \times \E_Y \times \E_Z @<\phi_1 << D @>\phi_2>> A @>\phi_3>> \E_V,  
\end{CD}
\tag{**}
\]
where we set $\phi_1:=q_3\, n_3, \phi_2:= n_2$ and $ \phi_3:= t$.

Note that $\phi_2$ is a $\G_X \times \G_Y \times \G_Z$-principal
bundle. In particular, the fiber dimension of $\phi_2$ is 
\[
f_2=\dim \G_X \times \G_Y \times \G_Z.
\]
The morphism $\phi_1$ is a smooth morphism with connected fibers of fiber
dimension
\begin{equation*}
\begin{split}
f_1&=\dim \G_X \times \G_Y \times \G_Z
+\sum_{i\in I} |Y_i| \, |Z_i| +|X_i|\, |Z_i| +|X_i|\, |Y_i|\\
& + \sum_{h\in H} |Y_{s(h)}|\, |Z_{t(h)}| +|X_{s(h)}|\, |Z_{t(h)}| +
|X_{s(h)}|\, |Y_{t(h)}|
\end{split}
\end{equation*}

Similarly, we can compute the fiber dimensions of the morphisms 
$p_1, p_2, q_1$ and $q_2$. They are 
\begin{enumerate}

\item[] $e_1=\dim \G_X \times \G_Y +\sum_{i\in I} |X_i|\, |Y_i|+
      \sum_{h\in H} |X_{s(h)}|\, |Y_{t(h)}|$;

\item[] $e_2=\dim \G_X \times \G_Y$;

\item[] $d_1=\dim \G_U \times \G_Z +\sum_{i\in I} |U_i|\, |Z_i|+
       \sum_{h\in H} |U_{s(h)}|\, |Z_{t(h)}|$; and 

\item[] $d_2=\dim \G_U \times \G_Z$, respectively.

\end{enumerate}
From the above analysis, one can check directly that 
\[
d=d_1-d_2+e_1-e_2=f_1-f_2.
\]
So we have 
\begin{cor}
$K \circ  L \circ  M = (\phi_3)_!\, (\phi_2)_{\flat}\, \phi_1^* 
(K \boxtimes L \boxtimes M)[f_1-f_2]$.
\end{cor}

More generally, consider the diagram
\[
\label{***}
\begin{CD}
\E_{V(1)} \times \cdots \times \E_{V(n)} @<\phi_1<< D_n @>\phi_2>> A_n @>\phi_3>> \E_V,
\end{CD}
\tag{***}
\]
where $V(m)$ ($m=1,\cdots, n$ and  $n\geq 2$) is an $I$-graded space such that
$V=\oplus_{m=1}^n \, V(m)$. 

$A_n$ is the variety consisting of all pairs $(x, \Vb)$, 
where $x\in \E_V$ and 
$\Vb=(V=V^0 \supseteq V^1 \supseteq \cdots \supseteq V^n=0)$,
such that $|V^{m-1}/V^m|= |V(m)|$ ($m=1,\cdots, n$) and 
$\Vb$ is $x$-stable.

$D_n$ is the variety consisting of all triples $(x, \Vb, \mg)$, 
where $(x,\Vb) \in A_n$ and 
$\mg$ is an $n$-tuple of linear isomorphisms 
$\mrm g^m: V^{m-1}/V^m \to V(m)$ ($ m=1,\cdots, n$).

The morphisms $\phi_1, \phi_2, \phi_3$ are defined similar to the morphisms in the
diagram (\ref{**}). Denote by $f_1^{(n)}, f_2^{(n)}$ the fiber
dimensions of $\phi_1,\phi_2$, respectively.

Assume that $K_m \in \mQ_{V(m)}$, for $m=1,\cdots n$.
We have 
\begin{cor}
\label{multn}
$K_1 \circ \cdots \circ K_n =  
(\phi_3)_!\, (\phi_2)_{\flat}\, (\phi_1)^* 
(K_1 \boxtimes \cdots \boxtimes K_n) \, [f_1^{(n)}-f_2^{(n)}]$.
\end{cor}

\begin{proof}
We prove by induction. 
when $n=2$, the statement in Corollary ~\ref{multn} holds
automatically. Assume that the statement in Corollary ~\ref{multn}
holds. Assume that $K_{n+1}\in \mQ_{V(n+1)}$. Let
$V=\oplus_{m=1}^{n+1}\, V(m)$, $W=V(n+1)$ and $T=\oplus_{m=1}^n\, V(m)$.
Then 
\begin{equation}
\label{I}
\begin{split}
&(K_1\circ \cdots \circ K_n) \circ K_{n+1}\\
&=(q_3)_!\, (q_2)_{\flat}\, (q_1)^*\, ((K_1\circ \cdots \circ K_n)
\boxtimes K_{n+1})[d_1-d_2]\\
&=(q_3)_!\, (q_2)_{\flat}\, (q_1)^*\, 
\left ((\phi_3)_!\, (\phi_2)_{\flat}\, (\phi_1)^* 
(K_1 \boxtimes \cdots \boxtimes K_n) 
\boxtimes K_{n+1}\right )[d]
\end{split}
\end{equation}
where $q_1, q_2$ and $q_3$ are the morphisms in the diagram (\ref{ind}),
$d_1$ and $d_2$ are the fiber dimensions of $q_1$ and $q_2$,
respectively and $d=d_1-d_2+f_1^{(n)}-f_2^{(n)}$.
Let $\text{id}$ be the identity morphism of $\E_{V(n+1)}$. 
The equation (\ref{I}) becomes
\begin{equation}
\label{II}
\begin{split}
&(K_1\circ \cdots \circ K_n) \circ K_{n+1}\\
&=(q_3)_!\, (q_2)_{\flat}\, (q_1)^*\, 
((\phi_3)_!\, (\phi_2)_{\flat}\, (\phi_1)^* 
(K_1 \boxtimes \cdots \boxtimes K_n) 
\boxtimes K_{n+1})[d])\\
&=(q_3)_!\, (q_2)_{\flat}\, (q_1)^*\, 
(\phi_3\times \text{id})_!\, (\phi_2\times \text{id})_{\flat}\,
(\phi_1\times \text{id})^* 
(\mathcal M)[d]\\
\end{split}
\end{equation}
where $\mathcal M=K_1 \boxtimes \cdots \boxtimes K_n 
\boxtimes K_{n+1}$.
Consider the diagram similar to diagram (\ref{***})
\[
\begin{CD}
\E_{V(1)} \times \cdots \times \E_{V(n+1)} @<\psi_1<< D_{n+1} @>\psi_2>> A_{n+1} @>\psi_3>> \E_V,
\end{CD}
\]
where $D_{n+1}$, $A_{n+1}$, $\psi_1, \psi_2$ and $\psi_3$ are defined
similar to the ones in diagram (\ref{***}).
Denote by $f_1^{(n+1)}$ and $f_2^{(n+1)}$ for the fiber dimensions of $\phi_1$ and $\phi_2$, respectively.
A procedure similar to the Proof of Lemma ~\ref{associativity} then
implies that 
\begin{equation}
\label{III}
\begin{split}
&(q_3)_!\, (q_2)_{\flat}\, (q_1)^*\, 
(\phi_3\times \text{id})_!\, (\phi_2\times \text{id})_{\flat}\,
(\phi_1\times \text{id})^* 
(\mathcal M)[d]\\
&=(\psi_3)_!\, (\psi_2)_!\, (\psi_1)^* \, (\mathcal M)[d].
\end{split}
\end{equation}
From (\ref{II}) and (\ref{III}), we have
\begin{equation*}
\begin{split}
K_1\circ \cdots \circ K_n \circ K_{n+1}
=(\psi_3)_!\, (\psi_2)_!\, (\psi_1)^* \, (\mathcal M)[d].
\end{split}
\end{equation*}
Finally, one can check that $d=d_1-d_2+f_1^{(n)}-f_2^{(n)}=f_1^{(n+1)}-f_2^{(n+1)}$.
Then the Corollary follows by induction.
\end{proof}

\subsection{Canonical basis}
\label{algebra}

Let $\mathcal{K}_V=\mathcal K(\mathcal Q_V)$ be the Grothendieck group of the category
$\mathcal Q_V$, i.e., it is   
the abelian group with one generator $\langle L\rangle$ for each isomorphism class
of objects in $\mathcal{Q}_V$ with relations: 
$\langle L\rangle+\langle L'\rangle=\langle L''\rangle$ 
if $ L'' \cong L\oplus L'$.

Let $v$ be an indeterminate. 
We set $\mathbb A=\mathbb Z[v,v^{-1}]$.
Define an $\mathbb A$-module structure on $\mathcal K_V$ by 
$v^n\langle L\rangle =\langle L[n]\rangle $ 
for any generator $\langle L\rangle\in\mathcal{Q}_V$ and 
$n\in \mathbb Z$.
From the construction, 
it is a free $\mathbb A$-module with basis $\langle L\rangle$
where $\langle L\rangle$ runs over $\mathcal P_V$.

From the construction, we have 
$\mathcal{K}_V \cong \mathcal{K}_{V'}$,
for any $V$ and $V'$ such that $|V|=|V'|$. 
For each $\nu\in \mathbb N[I]$, fix an $I$-graded vector space $V$ of dimension $\nu$.
Let 
\[
\mathcal{K}_{\nu}=\mathcal{K}_V,\quad 
\mathcal{K}=\oplus_{\nu\in \mathbb{N}[I]}\mathcal{K}_{\nu}
\quad \text{and} \quad
\mathcal{K}_Q=\mathbb{Q}(v)\otimes_{\mathbb A}\mathcal{K}.
\]
Also let
\[
\mathcal P_{\nu}=\mathcal P_V
\quad \text{and} \quad
\mathcal P=\cup_{\nu\in \mathbb N[I]} \mathcal P_{\nu}.
\]
For any $\alpha, \beta\in \mathbb N[I]$, the operation $\circ$ induces an 
$\mathbb A$-linear map
\[
\circ: \mathcal K_{\alpha} \otimes_{\mathbb A} \mathcal K_{\beta} 
\to \mathcal K_{\alpha +\beta}.
\]
By adding up these linear maps, we have a linear map
\[
\circ: \mathcal K \otimes_{\mathbb A} \mathcal K \to \mathcal K.
\] 
Similarly, the operation $\circ$ induces a $\mathbb Q(v)$-linear map
\[
\circ:\mathcal K_Q \otimes_{\mathbb Q(v)} \mathcal K_Q \to \mathcal K_Q.
\]
\begin{prop}
\label{algebraprop}
\begin{enumerate}
\item $(\mathcal K, \circ)$ $($resp. $(\mathcal K_Q, \circ))$ is an associative algebra 
      over $\mathbb A$ $($resp. $\mathbb Q(v))$.
\item $\mathcal P$ is an $\mathbb A$-basis of $(\mathcal K, \circ)$ 
      and a $\mathbb Q(v)$-basis of $(\mathcal K_Q, \circ)$. 
\end{enumerate}
\end{prop}
\begin{proof}
The associativity of $\circ$ follows from Lemma ~\ref{associativity}. 
\end{proof}

From now on, we simply write $\mathcal K$ (resp. $\mathcal K_Q$) for the algebra $(\mathcal K, \circ)$
(resp. $(\mathcal K_Q, \circ)$).
For any $m\leq n\in \mathbb N$, let 
\[
[n]=\frac{v^n-v^{-n}}{v-v^{-1}}, \quad 
[n]^!=\prod_{m=1}^n [m]
\quad
\text{and}
\quad
\bin{n}{m}=\frac{[n]^!}{[m]^![n-m]^!}.
\]

Let $\Gamma$ be the underlying graph of $Q$ and $(,)$ the symmetric Euler form defined in Section ~\ref{graph}.
Let $c_{ij}=(i,j)$ for any $i,j \in I$. 
Then $\mrm C=(c_{ij})_{i, j\in I}$ is a symmetric generalized Cartan matrix. 

Let $\tilde{\mbf U}^-$ be the free algebra over $\mathbb Q(v)$ generated by $F_i$ for all $i\in I$.
We set $F^{(n)}_i=\frac{F_i^n}{[n]^!}$ for all $i\in I$ and $n\in \mathbb N$.
Then the negative part $\mbf U^-$ of the quantized
enveloping algebra attached to $\mrm{C}$ is the quotient of $\tilde{\mbf U}^-$ 
by the two-sided ideal generated by 
\begin{equation*}
\sum_{p=0}^{1-c_{ij}} (-1)^p F_i^{(p)} F_j F_i^{(1-c_{ij}-p)}
\quad \text{for}\; i\neq j\in I.
\end{equation*}
Let $_{\mathbb A}\mbf U^-$ be the $\mathbb A$-subalgebra of $\mbf U^-$
generated by $F^{(n)}_i$ for $i\in I$ and $n\in \mathbb N$.

\begin{thm} $($~\cite{lusztig1}, ~\cite{lusztig2}, ~\cite{lusztig5}$)$
\label{iso}
The map $F_i^{(n)} \mapsto F(ni)$ for all $i\in I$ and $n\in\mathbb N$ induces an $\mathbb A$-algebra isomorphism
\[
_{\mathbb A}\phi: \,_{\mathbb A} \mbf U^- \to \mathcal K
\]
and a $\mathbb Q(v)$-algebra isomorphism
\[
\phi: \mbf U^- \to \mathcal K_Q.
\]
\end{thm}

By Theorem ~\ref{iso} and Proposition ~\ref{algebraprop}, 
the image $\mbf B$ of the set $\mathcal P$ under the map $_{\mathbb A}\phi$ (resp. $\phi$) is an $\A$-basis (resp. $\mathbb Q(v)$-basis) of 
$_{\A}\mbf U^-$ (resp. $\mathbf U^-$).
We call $\mbf B$ the {\bf canonical basis} of $_{\A} \mbf U^-$ (resp. $\mbf U^-$).
Under the isomorphism $_{\A} \phi$ (resp. $\phi$), we identify $_{\A}\mbf U^-$ (resp. $\mbf U^-$, $\mbf B$) with $\mathcal K$ 
(resp. $\mathcal K_Q$, $\mathcal P$). 

{\bf Remark.} The construction in this Section works for arbitrary quivers. See ~\cite{lusztig5} for a more general treatment that
works for any symmetrisable generalized Cartan matrix $\mrm C$. For a list of properties of the canonical basis $\mbf B$, 
see ~\cite{lusztig5}.

\section{Monomial basis in $\mbf U^-$}
\label{momomial}

In this Section, we assume that $Q$ has no oriented cycle. 
We study monomials in $_\A\mbf U^-=\mathcal K$ in some special order. 
These monomials will then produce a monomial basis for $\mbf U^-$ (in fact an $\A$-basis in $_\A\mbf U^-$, see Prop. ~\ref{C}), 
which will be further shown to coincide with  the monomial basis constructed in ~\cite{LXZ}.

\subsection{Preprojective and preinjective component}
\label{case1}
Recall that we order the vertex set $I$: $i_0, i_1, \cdots, i_n$ such that $i_j$ is a source of 
$\sigma_{i_{j+1}}\cdots \sigma_{i_{n-1}} \sigma_{i_n} Q$, $\cdots$.

For any $\alpha=\sum_{j=0}^n \alpha_{i_j}\ i_j \in \mathbb N[I]$, we set
\[
F(\alpha)=F(\alpha_{i_n}i_n)\circ \cdots \circ F(\alpha_{i_1} i_1)\circ F(\alpha_{i_0} i_0).
\]
Note that $F(\alpha_{i_j} i_j)$ has been defined in Example
~\ref{example}.
For convenience, we set $F(\alpha)=1$ if $\alpha=0$.
By Lemma ~\ref{ind1}, 
\[
F(\alpha)=L_{\s}=(\pi_{\s})_!(\bar{\Q}_l)[\dim \tF_{\s}]
\] 
where $\s=(\alpha_{i_n}i_n, \cdots, \alpha_{i_1}i_1, \alpha_{i_0}i_0)$.
By definition, $\mF_{\s}$ has only one element: 
\[
\mbf V^{\bullet}=(V\supseteq \oplus_{j=0}^{n-1} V_{i_j}\supseteq \cdots \supseteq V_{i_0}\supseteq 0).\]
Moreover, $\mbf V^{\bullet}$ is stabilized by every element in $\E_V$. 
So $\tF_{\s}\simeq \E_V$. Therefore 
\begin{lem}
\label{monomial1}
$F(\alpha)$ is a simple perverse sheaf.
\end{lem}

Denote by $rP_m$ the direct sum of $r$ copies of $P_m$ defined in Section ~\ref{classification}.
Then $rP_m$ has no self extension by Lemma ~\ref{vanishing}. 
In other words, the $\G_V$-orbit $O_{rP_m}$ corresponding to $rP_m$ is open in $\E_V$.
Thus
\begin{align}
\label{preprojective}
F(r|P_m|)=\mrm{IC}(O_{rP_m},\bar{\Q}_l).
\end{align}
Similarly
\begin{align}
\label{preinjective}
F(r|I_m|)=\mrm{IC}(O_{rI_m},\bar{\Q}_l).
\end{align}

\subsection{Inhomogeneous component}
\label{case2}
Let $\mathscr T$ be an inhomogeneous tube. Assume that $\mbf V=(V,x) \in \mrm{Rep}(\mathscr T)$ is aperiodic 
(see Section ~\ref{classification}). Denote by $O=O_x$ the $\G_V$-orbit of $x$ in $\E_V$. 
Let $\mrm{IC}(O,\bar{\Q}_l)$ be the simple perverse sheaf on $\E_V$ 
whose support is $\overline O$, the closure of $O$, and whose restriction to $O$ is $\bar{\Q}_l$.
\begin{prop}
\label{inhomogeneous1}
There exists a monomial, denote by $F(O)$ or $F(\mbf V)$, such that 
\begin{align}
\label{inhomogeneous}
F(O)=\mrm{IC}(O,\bar{\Q}_l)\oplus A
\end{align}
where $\mrm{Supp}(A)\subseteq \overline O\backslash O$.
\end{prop}

In fact, Lusztig showed in ~\cite{lusztig3} (see also ~\cite[Lem. 3.8]{Li-Lin}) that there exists a sequence
$(\gamma_{i_1},\gamma_{i_2},\cdots,\gamma_{i_m})$ such that
\begin{align}
F(\gamma_{i_1})\circ F(\gamma_{i_2})\circ \cdots \circ F(\gamma_{i_m})=\mrm{IC} (O,\bar{\Q}_l) \oplus A
\end{align}
where $A$ is a semisimple complex and $\gamma_{i_j}$ is one of the dimensions of the simple regular representations in $\mathscr T$.
By ~\cite[Thm. 8.5]{DDX}, the sequence can be chosen such that
$\mrm{Supp}(A) \subseteq \overline O \backslash O$. 
Proposition ~\ref{inhomogeneous1} follows.

\subsection{Homogeneous component}
\label{case3}
Let $\mbf 1=(1,\cdots ,1)$ be a partition of a positive integer $m$. We set 
$\mbf 1 \delta=(\delta,\cdots, \delta)$ such that $\delta+\cdots + \delta=m\delta$.
Let
\[
F(\mbf 1\delta)=F(\delta)\circ \cdots \circ F(\delta).
\]
Note that $F(\delta)$ is defined in Section ~\ref{case1}.
From Lemma ~\ref{monomial1}, $F(\mbf 1\delta)$ is in $\mathcal K$.
Moreover, by Lemma ~\ref{ind1} (2), 
\[
F(\mbf 1\delta)=L_{\mbf 1\delta}=(\pi_{\mbf 1\delta})_!(\bar{\Q}_l)[\dim \tF_{\mbf 1\delta}]
\]
where $\tF_{\mbf 1\delta}$ and $\pi_{\mbf 1\delta}: \tF_{\mbf 1\delta} \to \E_V$
is defined similar to $\tF_{\s}$ and  $\pi_{\s}$, respectively, in Section ~\ref{flag}.

Let $X(m\delta)$ be the subvariety of $\E_V$ consisting of all elements $x$ such that
\[
(V,x) \simeq T_1 \oplus\cdots \oplus T_m
\]
where $T_1, \cdots, T_m$ are pairwise non isomorphic regular simple representations of $Q$.
$X(m\delta)$ is an irreducible smooth variety  and open in its closure
$\overline X(m\delta)$ (see ~\cite{lusztig3}).

Denote by $\tilde{X}(m\delta)$ the variety consisting of all sequences $\mbf x=(x, T_{i_1}, \cdots, T_{i_m})$ where
$x\in X(m\delta)$ and $T_{i_1},\cdots , T_{i_m}$ is a permutation of $T_1,\cdots ,T_m$.
Denote by $\tau_{\mbf 1\delta}$ the first projection $\tilde{X} (m\delta) \to X(m\delta)$.
Then $\tau_{\mbf 1\delta}$  is a principal $\mathfrak S_m$-covering where
$\mathfrak S_m$ is the symmetric group of $m$ letters. So $\mathfrak S_m$ acts on the fiber
$\tau_{\mbf 1\delta}^{-1}(x)$ of $x$ in $X(m\delta)$, i.e.,
\begin{align}
\label{action}
s.\mbf x=\left (x, T_{s(i_1)}, \cdots, T_{s(i_m)}\right)
\quad \text{for all} \; s \in \mathfrak S_m\; \text{ and}\; \mbf x\in X(m\delta).
\end{align}

Let $M^{\mbf 1}(x)$ be the vector space over $\bar{\Q}_l$ spanned by elements in the fiber 
$\tau_{\mbf 1\delta}^{-1}(x)$. 
The action (\ref{action}) then induces an $\mathfrak S_m$-module structure on 
$M^{\mbf 1}(x)$.
By ~\cite{James}, $M^{\mbf 1}(x)$ is nothing but the permutation module of $\mathfrak S_m$ 
with respect to the partition $\mbf 1$. 
Denote by $\chi_{\boldsymbol{\mu}}$ the irreducible representation of $\mathfrak S_m$ 
corresponding to the partition $\boldsymbol{\mu}$ of $m$.
Then from ~\cite{James},
\begin{align}
\label{partition}
M^{\mbf 1}(x) = \chi_{\mbf 1} 
\oplus \oplus_{\boldsymbol{\mu}} A_{\mbf 1\boldsymbol{\mu}} \chi_{\boldsymbol{\mu}}
\end{align}
where $A_{\mbf 1\boldsymbol{\mu}}\chi_{\boldsymbol{\mu}}$ is a direct sum of $A_{\mbf 1\boldsymbol{\mu}}$ copies of
the representation $\chi_{\boldsymbol{\mu}}$ and $\boldsymbol{\mu}$ runs through all partitions of $m$ such that
$\boldsymbol{\mu}>\mbf 1$ and $>$ is the lexicography order for  the set of partitions of $m$ as follows.
Assume that $\boldsymbol{\lambda}=(\lambda_1,\cdots,\lambda_l)$ and $\boldsymbol{\mu}=(\mu_1,\cdots,\mu_{l'})$
are partitions of $m$, then
$\boldsymbol{\mu}>\boldsymbol{\lambda}$ if and only if there exist a $p$ such that $\mu_{p'}=\lambda_{p'}$ for all
$p'<p$ and $\mu_p>\lambda_p$.

Denote by $\mathcal L_{\boldsymbol{\lambda}}$ the local system on $X(\boldsymbol{\lambda}\delta)$ corresponding to
the irreducible representation $\chi_{\boldsymbol{\lambda}}$.
Since  $\tau_{\mbf 1\delta}$ is a principal $\mathfrak S_m$-covering, the monodromy representation of 
$(\tau_{\mbf 1\delta})_! (\bar{\Q}_l)$ at $x$ is $M^{\mbf 1}(x)$. 
From ~\cite{Iversen} and (\ref{partition}), 
\begin{align}
\label{local}
(\tau_{\mbf 1\delta})_!(\bar{\Q}_l)=\mathcal L_{\mbf 1} \oplus\oplus_{\boldsymbol{\mu}}
A_{\mbf 1\boldsymbol{\mu}}\mathcal L_{\boldsymbol{\mu}}
\end{align}
where the notations are defined similar to the notations in (\ref{partition}).

Given any $x\in X(m\delta)$, we have $(V,x)\simeq T_1\oplus \cdots \oplus T_m$. 
So we can decompose $V=\oplus_{r=1}^m V(r)$ such that $V(r)$ is $x$-stable and 
$|V(r)|=\delta$. We denote by $(V(r), x)$ the subreprsentation of $(V,x)$ by restricting $x$ to $V(r)$. 
By the fact that $|V(r)|=\delta$ and $(V(r), x)$ is a subrepresentation of a regular representation, we know that
$(V(r),x)\simeq T_{i_r}$ for some $i_r$. So we have $(V,x)\simeq \oplus_{r=1}^m (V(r),x)$. 
Moreover, this decomposable is unique up to order (see ~\cite[5.7]{Li-Lin}). 
Define an injective map $a: \tilde{X}(m\delta) \to \tF_{\mbf 1 \delta}$ by 
\[(x, T_1,\cdots, T_m) \mapsto (x, \mbf V^{\bullet})\]
where $\mbf V^{\bullet}=(V=V^0\supseteq V^1 \supseteq \cdots \supseteq V^m=0)$ with
$V^r=\oplus_{u=r}^m V(u)$ and $(V(u), x) \simeq T_u$ for $1\leq u\leq m$.
We have the following commutative diagram
\[
\begin{CD}
~\tilde{X}(m\delta) @>a>> \tF_{\mbf 1\delta}\\
@V\tau_{\mbf 1\delta}VV @V\pi_{\mbf 1\delta}VV\\
X(m\delta) @>b>> \E_V
\end{CD}
\]
where the bottom arrow is the nature injection. By construction, this diagram is Cartesian. So we have
\begin{align}
\label{m-delta}
b^* (\pi_{\mbf 1\delta})_!(\bar{\Q}_l)= (\tau_{\mbf 1 \delta})_! a^*(\bar{\Q}_l)=(\tau_{\mbf 1\delta})_!(\bar{\Q}_l).
\end{align}
Note that the support of $(\pi_{\mbf 1\delta})_! (\bar{\Q}_l)$ is  $\overline X(m\delta)$ and $X(m\delta)$ is open dense in
$\overline X(m\delta)$. 
Let $\mrm{IC}(X(m\delta), \mathcal L_{\boldsymbol{\mu}})$ be the simple
perverse sheaf on $\E_V$ whose support is $\overline X(m\delta)$ and whose
restriction to $X(m\delta)$ is $\mathcal L_{\boldsymbol{\mu}}[\dim X(m\delta)]$.

From (\ref{local}) and (\ref{m-delta}), we have

\begin{lem}
\label{mbf1}
$F(\mbf 1\delta)=\mrm{IC} (X(m\delta), \mathcal L_{\mbf 1}) \oplus \oplus_{\boldsymbol{\mu}}
A_{\mbf 1 \boldsymbol{\mu}} \mrm{IC}(X(m\delta), \mathcal L_{\boldsymbol{\mu}}) \oplus B$
where $\mrm{Supp}(B)$ $\subseteq$ $ \overline X(m\delta) \backslash X(m\delta)$.
\end{lem}

More generally, let $\boldsymbol{\lambda}=(\lambda_1,\cdots \lambda_l)$ be a partition of $m$.
Set $\boldsymbol{\lambda}\delta=(\lambda_1 \delta, \cdots, \lambda_l \delta)$.
Let
\begin{align}
\label{homogeneous}
F(\boldsymbol{\lambda}\delta)=F(\lambda_1\delta)\circ \cdots \circ F(\lambda_l\delta).
\end{align}
Similar to $F(\mbf 1\delta)$,
\[
F(\boldsymbol{\lambda}\delta)=L_{\boldsymbol{\lambda} \delta}=(\pi_{\boldsymbol{\lambda}\delta})_!(\bar{\Q}_l)[\dim \tF_{\boldsymbol{\lambda}\delta}].
\]
Denote by $\tilde{X}_{\boldsymbol{\lambda}}(m\delta)$ 
the variety consisting of all sequences
\[
\mbf x=(x, \oplus_{r=1}^{\lambda_1}T_{i_r}, \oplus_{r=\lambda_1+1}^{\lambda_1+\lambda_2} T_{i_r},
\cdots, \oplus_{r=\lambda_1+\cdots+\lambda_{l-1}+1}^mT_{i_r})
\] 
where $x\in X(m\delta)$ and $T_{i_1},\cdots , T_{i_m}$ is a permutation of $T_1,\cdots ,T_m$.
Denote by $\tau_{\boldsymbol{\lambda}\delta}$ the first projection $\tilde{X}_{\boldsymbol{\lambda}} (m\delta) \to X(m\delta)$.
Again $\mathfrak S_m$ acts on the fiber
$\tau_{\boldsymbol{\lambda}\delta}^{-1}(x)$ of $x$ in $X(m\delta)$, i.e.,
\begin{align}
\label{action2}
s.\mbf x=\left (x,\oplus_{r=1}^{\lambda_1}T_{s(i_r)}, \oplus_{r=\lambda_1+1}^{\lambda_1+\lambda_2} T_{s(i_r)},
\cdots, \oplus_{r=\lambda_1+\cdots+\lambda_{l-1}+1}^mT_{s(i_r)} \right)
\end{align}
for all  $ s \in \mathfrak S_m$ and $ \mbf x\in X(m\delta)$.

Let $M^{\boldsymbol{\lambda}}(x)$ be the vector space over $\bar{\Q}_l$ 
spanned by elements in the fiber  $\tau_{\boldsymbol{\lambda}\delta}^{-1}(x)$. 
The action (\ref{action2}) then induces an $\mathfrak S_m$-module structure on 
$M^{\boldsymbol{\lambda}}(x)$.
By ~\cite{James}, $M^{\boldsymbol{\lambda}}(x)$ is nothing but the permutation module of $\mathfrak S_m$ 
with respect to the partition $\boldsymbol{\lambda}$. 
Similar to (\ref{partition}), (\ref{local}) and (\ref{m-delta}), we have
\begin{align}
\label{lambda1}
&M^{\boldsymbol{\lambda}}(x) = \chi_{\boldsymbol{\lambda}} 
\oplus \oplus_{\boldsymbol{\mu}>\boldsymbol{\lambda}} A_{\boldsymbol{\lambda\mu}} \chi_{\boldsymbol{\mu}},\\
\label{lambda2}
&(\tau_{\boldsymbol{\lambda}\delta})_!(\bar{\Q}_l)=\mathcal L_{\boldsymbol{\lambda}} \oplus\oplus_{\boldsymbol{\mu}>\boldsymbol{\lambda}}
A_{\boldsymbol{\lambda\mu}}\mathcal L_{\boldsymbol{\mu}}, \quad \text{and}\\
\label{lambda3}
&b^*(\pi_{\boldsymbol{\lambda}\delta})_!(\bar{\Q}_l)=(\tau_{\boldsymbol{\lambda}\delta})_!(\bar{\Q}_l).
\end{align}
(Here $A_{\boldsymbol{\lambda\mu}}$ is a Kostka number.) From (\ref{lambda2}) and (\ref{lambda3}), we have
\begin{lem}
\label{lambda}
$F(\boldsymbol{\lambda}\delta)=\mrm{IC}(X(m\delta), \mathcal L_{\boldsymbol{\lambda}}) \oplus \oplus_{\boldsymbol{\mu}>\boldsymbol{\lambda}}
A_{\boldsymbol{\lambda} \boldsymbol{\mu}} 
\mrm{IC}(X(m\delta), \mathcal L_{\boldsymbol{\mu}}) \oplus B$
where $\mrm{Supp}(B)\subseteq \overline X(m\delta) \backslash X(m\delta)$.
\end{lem}

\subsection{Monomial basis}
\label{basis}
We construct a monomial basis for $\mbf U^-$ and show that it coincides with the one given in ~\cite{LXZ}.

For convenience, we say the sequence $(0)$ is a partition of $0$.
For any $\nu\in \mathbb N[I]$, let $\Delta_{\nu}$ be the set of all pairs 
$(\mbf a, \boldsymbol{\lambda})$, where $\mbf a: \Ind (Q) \to \mathbb N$ is a function and 
$\boldsymbol{\lambda}$ is a partition of some nonnegative integer $m$, 
satisfying the following conditions:
\begin{align}
\label{condition}
\sum_{\mbf V\in \Ind (Q)} \mbf a(\mbf V) |\mbf V| +m\delta=\nu
\quad \text{and}\quad
\prod_{r=0}^{p-1} \mbf a( (\Phi^+)^r (T))=0  
\end{align}
for any regular indecomposable representation $T$ of period $ p$. 

We set $\Delta=\sqcup_{\nu\in \mathbb N[I]} \Delta_{\nu}$.

Note that $\mbf a (T)=0$ for any regular indecomposable representation of period $1$.
For each inhomogeneous tube $\mathscr T_i$ ($i=1,\cdots s$), we set 
\[
\mbf a ^i=\oplus_{a, l} \, \mbf a (T_{i, a, l}) \, T_{i,a,l}
\]
where $a=1,\cdots, p_i$, $l\in \mathbb Z_{\geq 0}$ and $T_{i,a,l}$ are indecomposable in $\mathscr T_i$.
By (\ref{condition}), $\mbf a^i$ is aperiodic.

For simplicity, we write 
\[
\alpha_m=\mbf a(P_m)|P_m| 
\quad \text{and}\quad
\beta_l=\mbf a(I_l)|I_l|
\quad 
\text{for}\;  m\in \mathbb Z_{\geq 0}\;\text{and}\; l\in \mathbb Z_{\leq n}.
\]

We form monomials
\begin{align}
\label{pm}
&F(\mbf a^+)=F(\alpha_0)\circ F(\alpha_1)\circ \cdots \circ
F(\alpha_m)\circ \cdots\\
\label{pmm}
&F(\mbf a^-)=\cdots \circ F(\beta_l)\circ \cdots \circ F(\beta_{n-1})\circ F(\beta_n)
\end{align}
where $F(\alpha_m)$ and $ F(\beta_l)$ are defined 
in Section ~\ref{case1} (\ref{preprojective}) and  (\ref{preinjective}), respectively. 
The products are finite since $\mbf a$ is support finite by (\ref{condition}).
We set
\begin{align}
\label{monomial}
F(\mbf a, \boldsymbol{\lambda})
=F(\mbf a^+) \circ F(\mbf a^1)\circ \cdots \circ F(\mbf a^s)\circ 
F(\boldsymbol{\lambda}\delta) \circ F(\mbf a^-)
\end{align}
where $F(\mbf a^s)$ is defined in Section  ~\ref{case2} (\ref{inhomogeneous}) 
and $F(\boldsymbol{\lambda}\delta)$ is defined in Section ~\ref{case3} (\ref{homogeneous}).
 
\begin{prop}
\label{monomialbasis}
The monomials $F(\mbf a, \boldsymbol{\lambda})$ 
where $(\mbf a, \boldsymbol{\lambda}) \in \Delta$  
form a $\mathbb Q(v)$-basis for $\mathcal K_Q=\mbf U^-$.
Moreover, $\{F(\mbf a,\boldsymbol{\lambda})\;|\; (\mbf a,\boldsymbol{\lambda})\in \Delta\}$
can be identified with the monomial basis 
$\{\mbf{m_a}\; |\; \mbf a\in \mathcal M\}$ defined in ~\cite[Prop. 8.4]{LXZ}.
\end{prop}

\begin{proof}
We only show that  the set $\{F(\mbf a, \boldsymbol{\lambda})\;|\; (\mbf a,\boldsymbol{\lambda}) \in \Delta\}$ coincides with the monomial basis
$\{\mbf{m_c}\;|\: \mbf c\in \mathcal M\}$ defined in ~\cite[Prop. 8.4]{LXZ}, as monomials in $\mbf U^-$. 
Then the Proposition ~\ref{monomialbasis} follows from ~\cite[Prop. 7.6, 8.4]{LXZ}. 

Note that from ~\cite{LXZ}, $\mbf{m_c}$ can be written as 
\begin{align}
\label{monomial-LXZ}
\mbf{m_c}=\mbf{m_a}\star \mbf m_{\pi_1\mbf c}\star\cdots \star \mbf m_{\pi_s\mbf c} \star \mbf m_{\omega_{\mbf c}\delta}\star \mbf {m_b}.
\end{align}
Among the components in (\ref{monomial}) and (\ref{monomial-LXZ}), the identifications
\begin{align*}
F(\mbf a^+) \leftrightarrow \mbf{m_a},\quad
F(\mbf a^-)\leftrightarrow \mbf{m_b}\quad \text{and}\quad
F(\mbf a^i)\leftrightarrow \mbf{m}_{\pi_i\mbf c}
\end{align*}
are clear by comparing Section ~\ref{case1} in this paper with Section 6.4 in ~\cite{LXZ}, 
and Section ~\ref{case2} in this paper with Section 5.3 in ~\cite{LXZ}.
Now the only part needs to be identified is 
\begin{align}
\label{homogeneous-component}
F(\boldsymbol{\lambda}\delta) \leftrightarrow \mbf m_{\omega_{\mbf c}\delta}.
\end{align}
Recall from ~\cite{LXZ}, 
$\mbf m_{\omega_{\mbf c}\delta}=\mbf m_{\omega_1\delta} \star \mbf m_{\omega_2\delta} \star \cdots\star \mbf m_{\omega_l\delta}$
where $\omega_{\mbf c}=(\omega_1,\cdots,\omega_l)$ is a partition of some nonnegative integer. 
From (\ref{homogeneous}), we only need to identify 
\begin{align}
\label{homogeneous-LXZ}
F(\lambda\delta) \leftrightarrow \mbf m_{\lambda\delta}
\end{align}
for some positive integer $\lambda$.

Recall the definition  of $\mbf m_{\lambda \delta}$ from ~\cite{LXZ}. In fact, our definition is not exactly the same as the one
given in ~\cite{LXZ}, instead our definition gives a realization of the process of producing $\mbf m_{\lambda \delta}$ given in 
~\cite{LXZ}. There are four cases to consider.

{\bf Case 1.} 
When $Q$ is the Kronecker quiver $\bullet \rightrightarrows \bullet$, denoted by $K$. 
We denote by $\mbf i$ and $\mbf j$ for the sink and the source, $a$ and $b$ for the two arrows, respectively.
Then $\mbf m_{\lambda\delta}=E_{\mbf j}^{(\lambda)} \star E_{\mbf i}^{(\lambda)}$ in ~\cite{LXZ}.
But $E_{\mbf j}^{(\lambda)}$, $E_{\mbf i}^{(\lambda)}$ and the (generic) $\star$-product can be identified with
$F(\lambda\mbf j), F(\lambda \mbf i)$ and the $\circ$-product in this paper, respectively. 
So  $F(\lambda\delta)$ and $\mbf m_{\lambda\delta}$ can be identified in $\mbf U^-$ when $Q$ is the Kronecker quiver.

{\bf Case 2.}
Assume now that $Q$ is a quiver of type $A_n^{(1)}$ ($n\geq 3$), but not a cyclic quiver. 
Recall that $i_0$ is a sink. We denote by $h_1$ and $h_2$ for the two arrows that terminate at $i_0$.
Let $\mbf V=(V,x)$ be the representation of $Q$ defined by 
\begin{align*}
&V_i=k \quad \text{if} \; k \neq i\quad \text{and}\quad   V_{i_0}=0;\\
&x_h=\mrm{id}_{V_{s(h)}}  \quad \text{if}\; h \neq h_1, h_2\quad \text{and}\quad    x_h=0  \quad \text{if}\; h=h_1, h_2. 
\end{align*}
Then $\mbf V$ is an indecomposable representation of dimension $\delta-i_0$ and has no self-extension.
($\mbf V$ has no self-extension is  due to the fact that $\delta-i_0$ is a positive root of finite type.)

Note that a representation of $K$ is a quadruple $(V_{\mbf i}, V_{\mbf j}, x_a, x_b)$.
We define a functor 
\[
\phi: \mrm{Rep}(K) \to \mrm{Rep}(Q)\quad (V_{\mbf i},V_{\mbf j},x_a,x_b)\mapsto (W,y)
\]
by
\begin{align*}
&W_i=V_{\mbf j} \quad \text{if}\; i\neq i_0 \quad \text{and}\quad W_{i_0}=V_{\mbf i};\\
&y_h=x_h \quad \text{if}\; h\neq h_1, h_2 \quad \text{and}\quad y_{h_1}=x_a\quad \text{and}\quad y_{h_2}=x_b. 
\end{align*}
The assignments extend to a functor. 
From the definition, we have 
\begin{align*}
\phi(S_{\mbf j})=\mbf V \quad \text{and} \quad 
\phi(S_{\mbf i})=S_{i_0}.
\end{align*}
This functor can be identified with the functor $F$ defined in ~\cite{LXZ}.
One can check that $\phi$ is an exact embedding, so it induces an injective algebra homomorphism
\begin{align*}
\phi: \mbf U^-(K)\to \mbf U^-,
\end{align*}
where $\mbf U^-(K)$ (resp. $\mbf U^-$) is the negative part of 
the quantized enveloping algebra associated to the generalized Cartan matrix
of the underlying graph of $K$ (resp. $Q$).
Moreover, since $\mbf V$ has no self-extension,
\begin{align*}
\phi(E_{\mbf j}^{(\lambda)}) = E_{i_n}^{(\lambda)}\star E_{i_{n-1}}^{(\lambda)}\star\cdots\star
E_{i_1}^{(\lambda)}.
\end{align*}
Also we have $\phi(E_{\mbf i}^{(\lambda)})=E_{i_0}^{(\lambda)}$. By definition in ~\cite{LXZ},
\begin{align*}
\mbf m_{\lambda \delta}=\phi(E_{\mbf j}^{(\lambda)}) \star \phi(E_{\mbf i}^{(\lambda)})
=E_{i_n}^{(\lambda)}\star \cdots \star E_{i_1}^{(\lambda)} \star E_{i_0}^{(\lambda)}.
\end{align*}
From this and Section ~\ref{case3} in this paper, we can identify $F(\lambda\delta)$ with $\mbf m_{\lambda\delta}$.
The identification of (\ref{homogeneous-LXZ}) holds in this case.

{\bf Case 3.}
Assume that $Q$ is of type $D_n^{(1)}$ ($n\geq 4$) or $E_n^{(1)}$ ($n=6,7,8$). We further assume that
there is an extending vertex $i$ in $I$ such that $i$ is a sink. Since $i$ is a sink, when we order the vertex set $I$,
we can set $i_0=i$. 
Since $Q$ is of type $D_{n}^{(1)}$ or $E_n^{(1)}$, we know that there is only one arrow $h_0$ such that $t(h_0)=i_0$.
Let $\mbf V=(V,x)$ be an indecomposable representation of dimension $\delta-i_0$
and then has no self-extension. 
Such a representation exists since $\delta-i_0$ is a positive root of finite type.
Note that $(\delta-i_0)_{s(h_0)}=2$, i.e., $V_{s(h_0)}=k^2$. We fix a basis $\{e_1,e_2\}$ for $V_{s(h_0)}$.

Define a functor 
\[
\phi: \mrm{Rep}(K)\to \mrm{Rep}(Q)\quad (V_{\mbf i}, V_{\mbf j}, x_a,x_b) \mapsto (W,y)
\]
by
\begin{align*}
&W_i=V_{\mbf j}\otimes V_i\quad \text{if} \; i\neq i_0 \quad \text{and}\quad W_{i_0}=V_{\mbf i};\\
&y_h=\text{id}_{V_{\mbf j}}\otimes x_h: V_{\mbf j}\otimes V_{s(h)} \to V_{\mbf j}\otimes V_{t(h)} 
\quad \text{if}\; h\neq h_0 \quad \text{and}\\
&y_{h_0}=x_a\oplus x_b: V_{\mbf j} \otimes V_{s(h_0)} \to V_{\mbf i}.
\end{align*}
Here $x_a \oplus x_b$ is defined by 
\begin{align*}
x_a\oplus x_b (v\otimes e_1+ w \otimes e_2)=x_a(v)+x_b(w) \quad \text{for any } v, w \in V_{\mbf j}. 
\end{align*}
From definition, we have
\[
\phi(S_{\mbf j})=\mbf V \quad \text{and} \quad \phi(S_{\mbf i})=S_{i_0}.
\]
As in Case 1 and 2, $\phi$ induces an injective algebra homomorphism
\begin{align*}
\phi: \mbf U^-(K)\to \mbf U^-
\end{align*}
such that 
\begin{align*}
\phi(E_{\mbf i}^{(\lambda)})=E_{i_0}^{(\lambda)}
\quad \text{and}\quad
\phi(E_{\mbf j}^{(\lambda)})=E_{i_n}^{(\lambda \delta_{i_n})} \star \cdots \star E_{i_1}^{(\lambda\delta_{i_1})}.
\end{align*}
By the definition in ~\cite{LXZ}, 
\begin{align*}
\mbf m_{\lambda\delta}=\phi(E_{\mbf j}^{(\lambda)}) \star \phi(E_{\mbf i}^{(\lambda)})=
E_{i_n}^{(\lambda \delta_{i_n})}\star \cdots \star E_{i_1}^{(\lambda\delta_{i_1})} \star E_{i_0}^{(\lambda\delta_{i_0})}.
\end{align*}
By comparing the definition of $F(\lambda\delta)$, we see that the identification in (\ref{homogeneous-LXZ}) holds in this case.

{\bf Case 4.} 
Finally, we assume that $Q$ is a quiver of type $D_n^{(1)}$ and $E_n^{(1)}$. 
Assume that all extending vertices in $Q$ are sources. 
We fix one of the extending vertices and set it equal to $i_n$. 
Note that in our ordering of $I$, $i_n$ is a source in $Q$.
Let $\mbf V=(V,x)$ be an indecomposable representation of  dimension $\delta-i_n$
and has no self-extension. Let $h_n$ be the arrow such that $s(h_n)=i_n$.
Define a functor 
\[
\phi: \mrm{Rep}(K)\to \mrm{Rep}(Q)\quad (V_{\mbf i}, V_{\mbf j}, x_a,x_b) \mapsto (W,y)
\]
by
\begin{align*}
&W_i=V_{\mbf i}\otimes V_i\quad \text{if} \; i\neq i_n \quad \text{and}\quad W_{i_n}=V_{\mbf j};\\
&y_h=\text{id}_{V_{\mbf i}} \quad \text{if}\; h\neq h_0 \quad \text{and}\quad 
y_{h_0}=[x_a, x_b]: V_{\mbf j}\to   V_{\mbf i}\otimes k^2.
\end{align*}
Here $[x_a, x_b]$ is defined by 
\begin{align*}
[x_a, x_b](v)= x_a(v)\otimes e_1+ x_b(v) \otimes e_2  \quad \text{for any } v \in V_{\mbf j}. 
\end{align*}
From the definition, we have
\[
\phi(S_{\mbf j})=S_{i_n} \quad \text{and} \quad \phi(S_{\mbf i})=\mbf V.
\]
As in Case 1, 2 and 3, $\phi$ induces an injective algebra homomorphism
\begin{align*}
\phi: \mbf U^-(K)\to \mbf U^-
\end{align*}
such that 
\begin{align*}
\phi(E_{\mbf i}^{(\lambda)})=E_{i_{n-1}}^{(\lambda \delta_{i_n})} \star \cdots \star E_{i_0}^{(\lambda\delta_{i_0})}
\quad \text{and}\quad
\phi(E_{\mbf j}^{(\lambda)})= E_{i_n}^{(\lambda)}.
\end{align*}
By the definition in ~\cite{LXZ}, 
\begin{align*}
\mbf m_{\lambda\delta}=\phi(E_{\mbf j}^{(\lambda)}) \star \phi(E_{\mbf i}^{(\lambda)})=
E_{i_n}^{(\lambda \delta_{i_n})}\star E_{i_{n-1}}^{(\lambda\delta_{i_{n-1}})} \star \cdots \star E_{i_0}^{(\lambda\delta_{i_0})}.
\end{align*}
Again from the definition of $F(\lambda\delta)$, we see that the identification in (\ref{homogeneous-LXZ}) holds in this case.
This finishes the Proof of Proposition ~\ref{monomialbasis}.
\end{proof}

{\bf Remark.} The functors $\phi$ constructed in Proposition ~\ref{monomialbasis} are in ~\cite{FMV}. 
The fact that $\{F(\mbf a,\boldsymbol{\lambda})\;|\; (\mbf a,\boldsymbol{\lambda})\in \Delta \}$ 
is a monomial basis of $\mbf U^-$ can be proved directly using the machinery 
built up in ~\cite{Li-Lin}, as we will see in the Proof of Theorem ~\ref{monomial-canonical}.

\section{Main results}
\label{main}

\subsection{}
\label{mainresults}
Let $V$ be an $I$-graded space over $k$ of dimension $\nu$. 
For any $(\mbf a,\boldsymbol{\lambda})\in \Delta_{\nu}$ (see ~\ref{basis}) 
with $\boldsymbol{\lambda}$ a partition of some nonnegative integer $m$, 
define $X(\mbf a,\boldsymbol{\lambda})$ to be the subvariety of $\E_V$
consisting of all elements $x$ in $\E_V$ such that
\[
(V,x)\simeq \oplus_{\mbf V\in \Ind (Q)} \mbf a(\mbf V) \mbf V 
\oplus T_1 \oplus\cdots \oplus T_m
\]
where $T_1,\cdots, T_m$  are pairwise non isomorphic simple regular representations of $Q$.
Then $X(\mbf a,\boldsymbol{\lambda})$ is a smooth irreducible variety 
and is open dense in its closure
$\overline X(\mbf a,\boldsymbol{\lambda})$ (see ~\cite{lusztig3}).

Define $\tilde  X(\mbf a,\boldsymbol{\lambda})$ to be the variety 
consisting of all sequences $\mbf x=(x,T_{r_1},\cdots,T_{r_m})$ where
$x\in \E_V$ and $T_{r_1},\cdots,T_{r_m}$ is a permutation of $T_1,\cdots T_m$.

Let $\pi_{\mbf a,\boldsymbol{\lambda}}$ be the first projection 
$\tilde X(\mbf a,\boldsymbol{\lambda})\to X(\mbf a,\boldsymbol{\lambda})$. 
$\pi_{\mbf a,\boldsymbol{\lambda}}$ is a principal $\mathfrak S_m$-covering (\cite{lusztig3}).
Thus the fundamental group of $X(\mbf a,\boldsymbol{\lambda})$ is $\mathfrak S_m$.
For each partition $\boldsymbol{\mu}$ of $m$, it determines 
an irreducible representation of $\mathfrak S_m$, denoted by 
$\chi_{\boldsymbol{\mu}}$. 
Let $\mathcal L_{\boldsymbol{\mu}}$ be 
the irreducible local system on $X(\mbf a,\boldsymbol{\lambda})$ such that 
the monodromy representation induced by $\mathcal
L_{\boldsymbol{\mu}}$ at the stalk of $x\in X(\mbf a,\boldsymbol{\lambda})$ is 
$\chi_{\boldsymbol{\mu}}$ (see ~\cite{Iversen}).

Denote by $\mrm{IC}(\mbf a,\boldsymbol{\lambda})$ the simple perverse
sheaf on $\E_V$ whose support is 
$\overline X(\mbf a,\boldsymbol{\lambda})$ and whose restriction to 
$X(\mbf a,\boldsymbol{\lambda})$
is $\mathcal L_{\boldsymbol{\lambda}}[\dim X(\mbf a,\boldsymbol{\lambda})]$.

We define a partial order, $\prec$, on the set $\Delta_{\nu}$ by
\begin{align*}
(\mbf a,\boldsymbol{\lambda}) \prec (\mbf b,\boldsymbol{\mu})
\quad \text{if and only if } & X(\mbf a,\boldsymbol{\lambda}) \subseteq 
 \overline X(\mbf b,\boldsymbol{\mu})\backslash X(\mbf b,\boldsymbol{\mu}) \quad \text{or}\\
&\overline X(\mbf a,\boldsymbol{\lambda})=\overline X(\mbf b,\boldsymbol{\mu}) \;\text{and}\; 
\boldsymbol{\mu}<\boldsymbol{\lambda}. 
\end{align*}
Here $<$ is the lexicographic order on the set of all partitions of $m$.
Note that the partial order is well-defined since if 
$\overline X(\mbf a,\boldsymbol{\lambda})=\overline X(\mbf b,\boldsymbol{\mu})$
then $\mbf a=\mbf b$.
The partial order $\prec$ on $\Delta_{\nu}$ then extends to a partial
order on $\Delta$, denoted again by $\prec$.
Now we can state  our main results in this paper.

\begin{thm}
\label{monomial-canonical}
Under the partial order $\prec$ on $\Delta$,
the transition matrix between the set 
$\mathcal F=
\{F(\mbf a,\boldsymbol{\lambda})\;|\;(\mbf a,\boldsymbol{\lambda}) \in \Delta\}$ and the set
$\mathcal{IC}=
\{\mrm{IC}(\mbf a,\boldsymbol{\lambda})\;|\;(\mbf a,\boldsymbol{\lambda})\in \Delta\}$  is 
upper triangular with entries in the diagonal equal to $1$ and 
entries in the upper diagonal in $\mathbb Z_{\geq 0}[v,v^{-1}]\subseteq \A$. 
More precisely,
\[
F(\mbf a, \boldsymbol{\lambda})= 
\mrm{IC}(\mbf a, \boldsymbol{\lambda}) + 
\sum_{\boldsymbol{\lambda} < \boldsymbol{\mu}}
A_{\boldsymbol{\lambda \mu}} \mrm{IC}(\mbf a, \boldsymbol{\mu}) +K,  
\]
where $A_{\boldsymbol{\lambda \mu}}$ is a Kostka number and 
$K$ is a linear combination 
(with coefficients in $\mathbb Z_{\geq  0}[v,v^{-1}]$)  
of elements in $\mathcal{IC}$
supported on 
$\overline X(\mbf a, \boldsymbol{\lambda})\backslash X(\mbf a,\boldsymbol{\lambda})$.
\end{thm}

The proof of Theorem ~\ref{monomial-canonical} will be given in Section ~\ref{Proof}.

\begin{prop}
\label{B}
$\mbf B =\mathcal{IC}.$
\end{prop}
Note that we identify $\mbf B$ with $\mathcal P$ in Section ~\ref{algebra}.
From Theorem ~\ref{monomial-canonical}, we have 
$\mrm{IC}(\mbf a,\boldsymbol{\lambda})\in \mathcal P$. 
Now by Proposition ~\ref{monomialbasis} and Theorem
~\ref{monomial-canonical}, Proposition ~\ref{B} follows.

As a consequence of Theorem ~\ref{monomial-canonical} and Proposition ~\ref{B}, we have
\begin{prop}
\label{C}
The basis $\mathcal F=\{F(\mbf a,\boldsymbol{\lambda})\;|\; (\mbf a,\boldsymbol{\lambda})\in \Delta\}$ is an $\A$-basis of 
$_\A\mbf U^-=\mathcal K$.
\end{prop}

Recall from ~\cite[8.4]{LXZ}, that the transition matrix between the monomial basis
$\{\mbf{m_c}\;|\; \mbf c\in \mathcal M\}$ ($=\mathcal F$) and the PBW-basis
$E_{\mathcal M} =\{E^{\mbf c}\;|\; \mbf c\in \mathcal M\}$ is upper triangular with entries in the diagonal equal to $1$ 
and entries above diagonal in $\mathbb Q[v,v^{-1}]$.
Since we can identify the index set $\mathcal M$ in ~\cite{LXZ} with $\Delta$ in this paper, we write 
$E^{(\mbf a,\boldsymbol{\lambda})}$ (resp. $\mbf m_{\mbf a,\boldsymbol{\lambda}}$) 
in this paper for $E^{\mbf c}$ (resp. $\mbf{m_c}$) in ~\cite{LXZ}.
(In fact $\mbf m_{\mbf a,\boldsymbol{\lambda}}=F(\mbf a,\boldsymbol{\lambda})$ by Proposition ~\ref{monomialbasis}.)
Combining this with Theorem ~\ref{monomial-canonical} and Proposition ~\ref{B}, we have

\begin{cor}
\label{D}
The transition matrix between the canonical basis $\mbf B$ and 
the PBW-basis  $\{E^{(\mbf a,\boldsymbol{\lambda})}\;|\; (\mbf a,\boldsymbol{\lambda})\in \Delta\}$
is upper triangular with entries in the diagonal equal to $1$ and
entries above the diagonal in $\mathbb Q[v,v^{-1}]$.
\end{cor}

\subsection{Proof of Theorem ~\ref{monomial-canonical}}
\label{Proof}
The subsection is devoted to proving Theorem ~\ref{monomial-canonical}.
The proof is essentially the same as the Proof of Proposition 5.10 in
~\cite{Li-Lin}, which comes from ~\cite{lusztig3}.

Let $V(\alpha_m), V(\beta_l), V(\mbf a^i)$ and $V(\boldsymbol{\lambda\delta})$ be the 
vector spaces of dimension $\alpha_m, \beta_l, |\mbf a^i|$ and $m\delta$, respectively,
where $m\in \mathbb Z_{\geq 0}, l\in \mathbb Z_{\leq n}, i=1,\cdots,
s$ and $\boldsymbol{\lambda}$ is a partition of $m$.
Denote by $V$ the direct sum of 
$V(\alpha_m), V(\beta_l), V(\mbf a^i)$ and $V(\boldsymbol{\lambda}\delta)$ for 
$m\in \mathbb Z_{\geq 0}, l\in \mathbb Z_{\leq n}$ and $i=1, \cdots, s$.
For simplicity, we write $\E_{\alpha_m}$, 
$\E_{\beta_l}$, $\E_{\mbf  a^i} $ 
and $\E_{\boldsymbol{\lambda}\delta}$ for the varieties
$\E_{V(\alpha_m)}$, $\E_{V(\beta_l)}$, $\E_{V(\mbf a^i)}$ and 
$\E_{V(\boldsymbol{\lambda}\delta)}$, respectively.
We set
\begin{align*}
\E_{\mbf a,\boldsymbol{\lambda}}
=\E_{\alpha_0}\times \cdots\times \E_{\alpha_m}\times \cdots \times
\E_{\mbf a^1}\times \cdots \times \E_{\mbf a^s} \times 
\E_{\boldsymbol{\lambda}\delta}\times \cdots \times
\E_{\beta_l}\times \cdots \times \E_{\beta_n}.
\end{align*}
Similarly, we define $\G_{\mbf a,\boldsymbol{\lambda}}$ 
for the product of the various general linear groups $\G_{V(m)}$.
Applying the diagram (\ref{***}) in Section ~\ref{induction}, we have the following diagram
\[
\label{d1}
\begin{CD}
\E_{\mbf a,\boldsymbol{\lambda}} @<\phi_1<< D_{\infty} @>\phi_2>> A_{\infty} @>\phi_3>> \E_V
\end{CD}
\tag{*}
\]
where $D_{\infty}$ and $A_{\infty}$ are defined similar to $D_n$ and $A_n$, respectively, in diagram (\ref{***})
in Section ~\ref{induction}. 
So are the maps $\phi_1, \phi_2$ and $\phi_3$. 
Denote by $f_1^{\infty}$ and $f_2^{\infty}$ 
for the fibre dimensions of $\phi_1$ and $\phi_2$, respectively.
By Lemma ~\ref{Reineke}, we have 
\[\phi_3(A_{\infty})=\overline X(\mbf a, \boldsymbol{\lambda}).\]
We set
\begin{align*}
&F^+=F(\alpha_0)\boxtimes F(\alpha_1)\boxtimes 
\cdots\boxtimes F(\alpha_m) \boxtimes \cdots,\\
&F^-= 
\cdots \boxtimes F(\beta_l)\boxtimes \cdots \boxtimes 
F(\beta_{n-1})\boxtimes F(\beta_n), \quad \text{and}\\
&F_{\mbf a,\boldsymbol{\lambda}}=F^+ \boxtimes
F(\mbf a^1)\boxtimes \cdots \boxtimes F(\mbf a^s) \boxtimes 
F(\boldsymbol{\lambda}\delta)\boxtimes F^-.
\end{align*}
By Corollary ~\ref{multn}, we have
\begin{align}
\label{eq1}
F(\mbf a,\boldsymbol{\lambda})=
(\phi_3)_! (\phi_2)_{\flat} (\phi_1)^* (F_{\mbf a,\boldsymbol{\lambda}})
[f_1^{\infty}-f_2^{\infty}].
\end{align}
Let $O_{\alpha_m}, O_{\beta_l}$ and $O_{\mbf a^i}$ be the $\G_{V(m)}$-orbits in 
$E_{\alpha_m}, E_{\beta_l}$ and $E_{\mbf a^i}$ corresponding
to  the representations $\mbf a(P_m)P_m$, $\mbf a(I_l)I_l$ and $\mbf a^i$, respectively. 
Let $O_{\boldsymbol{\lambda}\delta}$ be $X(m\delta)$ defined in Section ~\ref{case3}.
Denote by 
\begin{align*}
&O^+=O_{\alpha_0}\times O_{\alpha_1} \times \cdots \times O_{\alpha_m}\times \cdots,\\
&O^-=\cdots \times O_{\beta_l} \times \cdots \times 
O_{\beta_{n-1}} \times O_{\beta_n}, \quad \text{and}\\
&O_{\mbf a,\boldsymbol{\lambda}}=O^+\times O_{\mbf a^1}\times \cdots\times O_{\mbf a^s} \times 
O_{\boldsymbol{\lambda}\delta}  \times O^-.
\end{align*}
By (\ref{preprojective}), (\ref{preinjective}), (\ref{inhomogeneous}) and Lemma ~\ref{lambda}, 
$O_{\alpha_m}, O_{\beta_l}, O_{\mbf a^i}$ and 
$O_{\boldsymbol{\lambda}\delta}$ are open dense in
the supports of the complexes $F(\alpha_m), F(\beta_l), F(a^i)$ 
and $F(\boldsymbol{\lambda}\delta)$, respectively, where
$m\in \mbb Z_{\geq 0}, l\in \mbb Z_{\leq n}$ and $i=1, \cdots, s$.

Denote by $\mathcal N_{\alpha_m}, \mathcal N_{\beta_l}, \mathcal N_{\mbf a^i}$ 
and $\mathcal N_{\boldsymbol{\lambda}\delta}$
for the local system by restricting $F(\alpha_m)$, $F(\beta_l)$, 
$F(\mbf a^i)$ and $F(\boldsymbol{\lambda}\delta)$ to 
$O_{\alpha_m}, O_{\beta_l}, O_{\mbf a^i}$ and $O_{\boldsymbol{\lambda}\delta}$, respectively.
Note  that from (\ref{preprojective}), (\ref{preinjective}) and
(\ref{inhomogeneous}),  all local systems 
are trivial except $\mathcal N_{\boldsymbol{\lambda}\delta}$. 
By (\ref{lambda2}) and (\ref{lambda3}) in Section ~\ref{case3}, we have
\begin{align}
\label{eq2}
\mathcal N_{\boldsymbol{\lambda}\delta}=
\mathcal L_{\boldsymbol{\lambda}}\oplus \oplus_{\boldsymbol{\lambda}<\boldsymbol{\mu}}
A_{\boldsymbol{\lambda\mu} }\mathcal L_{\boldsymbol{\mu}}.
\end{align}

Let
$\mathcal N_{\mbf a,\boldsymbol{\lambda}}$ be the local system on 
$O_{\mbf a,\boldsymbol{\lambda}}$ by tensoring the local systems
$\mathcal N_{\alpha_m}$, $\mathcal N_{\beta_l}$, 
$\mathcal N_{\mbf a^i}$ and $\mathcal N_{\boldsymbol{\lambda}\delta}$,
i.e., $\mathcal N_{\mbf a,\boldsymbol{\lambda}}
= F_{\mbf a,\boldsymbol{\lambda}}|_{O_{\mbf a,\boldsymbol{\lambda}}}$,
the restriction of $F_{\mbf a,\boldsymbol{\lambda}}$ to $O_{\mbf a,\boldsymbol{\lambda}}$.

Consider the following diagram
\[
\label{d2}
\begin{CD}
O_{\mbf a,\boldsymbol{\lambda}} @<\phi_1'<<  D_{\infty}' @>\phi_2'>>
A_{\infty}' @>\phi_3'>> X(\mbf a,\boldsymbol{\lambda})\\
@VaVV @VbVV @VcVV @VdVV\\
\E_{\mbf a,\boldsymbol{\lambda}} @<\phi_1<< D_{\infty} @>\phi_2>>
A_{\infty} @>\phi_3>> \overline X(\mbf a,\boldsymbol{\lambda}) 
\end{CD}
\tag{**}
\]
where  $A_{\infty}'=\phi_3^{-1} (X(\mbf a,\boldsymbol{\lambda}))$,
$D_{\infty}'=(\phi_2)^{-1} (A_{\infty}')$, the bottom line is diagram
(\ref{d1}) in this subsection, the horizontal arrows in the top line
is induced from the arrows in the bottom line. 
The vertical maps are open embeddings.
Note that $\phi_1'$ is well-defined due to Lemma ~\ref{Reineke}. The
squares  on the right and in the middle are cartesians. From this and
diagram (\ref{d2}), one has
\begin{align}
\label{eq4}
\begin{split}
d^* F(\mbf a,\boldsymbol{\lambda}) 
&= d^* (\phi_3)_! (\phi_2)_{\flat} (\phi_1)^* (F_{\mbf a,\boldsymbol{\lambda}})
[f_1^{\infty}-f_2^{\infty}]\\
&=(\phi_3')_! (\phi_2')_{\flat} (\phi_1')^* 
(F_{\mbf a,\boldsymbol{\lambda}}|_{O_{\mbf a, \boldsymbol{\lambda}}})
[f_1^{\infty}-f_2^{\infty}].
\end{split}
\end{align}

Let $\tilde O_{\mbf a,\boldsymbol{\lambda}}$ be the variety defined 
as $\tilde X (m\delta)$  in Section ~\ref{case3}  with
$O_{\boldsymbol{\lambda}\delta}$ replaced by $O_{\mbf a,\boldsymbol{\lambda}}$.
Let $\pi$ be the first projection 
$\tilde  O_{\mbf a,\boldsymbol{\lambda}}\to O_{\mbf a,\boldsymbol{\lambda}}$.
Consider the following commutative diagram
\[
\label{d3}
\begin{CD}
\tilde O_{\mbf a,\boldsymbol{\lambda}} @<<< \tilde D_{\infty}' 
@>>> \tilde A_{\infty}' @>>> \tilde X(\mbf a,\boldsymbol{\lambda}) \\
@V\pi VV @VVV @VVV @V\pi_{\mbf a,\boldsymbol{\lambda}} VV \\
O_{\mbf a,\lambda} @<\phi_1'<< D_{\infty}' @>\phi_2'>> A_{\infty}' 
@>\phi_3'>> X(\mbf a, \boldsymbol{\lambda})
\end{CD}
\tag{***}
\]
where 
$\tilde A_{\infty}'=(\phi_3')^{-1}(\tilde X(\mbf a,\boldsymbol{\lambda}))$,
$\tilde D_{\infty}'=(\phi_2')^{-1}(\tilde A_{\infty}')$, and the maps
in the bottom are defined in diagram (\ref{d2}) and the other maps are
defined in the natural way. Note that the squares on the right and in the middle are
cartesians. From this and diagram (\ref{d3}), one has
\begin{align}
\label{eq5}
(\pi_{\mbf a,\boldsymbol{\lambda}})_!(\bar{\Q}_l)
=(\phi_3')_! (\phi_2')_{\flat} (\phi_1')^* ( \pi_!(\bar{\Q}_l))
=(\phi_3')_! (\phi_2')_{\flat} (\phi_1')^* 
(F_{\mbf a,\boldsymbol{\lambda}}|_{O_{\mbf a,\boldsymbol{\lambda}}}).
\end{align}
Combining (\ref{eq4}) and (\ref{eq5}), 
\begin{align*}
d^*F(\mbf a, \boldsymbol{\lambda})=
(\pi_{\mbf a,\boldsymbol{\lambda}})_!(\bar{\Q}_l)
[f_1^{\infty}-f_2^{\infty}].
\end{align*}
Note that $X(\mbf a,\boldsymbol{\lambda})$ is open dense in 
$\overline X(\mbf a, \boldsymbol{\lambda})$. The above equation
implies that up to shift,
\begin{align}
\label{eq6}
\begin{split}
F(\mbf a, \boldsymbol{\lambda})
&=\mrm{IC}( 
X(\mbf a, \boldsymbol{\lambda}), (\pi_{\mbf a,\boldsymbol{\lambda}})_!(\bar{\Q}_l))
\oplus K\\
&=\mrm{IC}(\mbf a, \boldsymbol{\lambda}) \oplus
\oplus_{\boldsymbol{\lambda} < \boldsymbol{\mu}}
A_{\boldsymbol{\lambda \mu}} \mrm{IC}(\mbf a, \boldsymbol{\mu}) \oplus K,
\end{split}
\end{align}
where $K$ is a semisimple complex supported on 
$\overline X(\mbf a, \boldsymbol{\lambda}) \backslash X(\mbf
a,\boldsymbol{\lambda})$.
But $F(\mbf a, \boldsymbol{\lambda})$ and 
$\mrm{IC}(\mbf a, \boldsymbol{\lambda})$ are stable under the bar
involution, the shift must disappear. Thus equation (\ref{eq6}) holds
without shift. 
(In other words, 
$\dim X(\mbf a,\boldsymbol{\lambda})= f_1^{\infty}-f_2^{\infty}$.)
This finishes the proof of Theorem ~\ref{monomial-canonical}.

\subsection{Resolutions of singularities}
\label{resolution}
In this subsection, we give a resolution of 
$\overline X(\mbf a, \boldsymbol{\lambda})$ for any 
$(\mbf a,\boldsymbol{\lambda}) \in \Delta$.
Note that $\overline X(\mbf a, \boldsymbol{\lambda})$ are the
supports for the affine canonical basis elements (as simple perverse
sheaves).

We preserve the setting in Section ~\ref{momomial}. For any
$\alpha=\sum_{j=0}^n \alpha_{i_j} i_j\in \mathbb N[I]$, we form a
sequence
\begin{align*}
\mbf s(\alpha)=(\alpha_{i_n} i_n, \cdots, \alpha_{i_1}i_1, \alpha_{i_0}i_0),
\end{align*}
where $i_0, i_1, \cdots, i_n$ is a fixed sequence of vertices in $I$
in Section ~\ref{momomial}. 

For any pair $(\mbf a, \boldsymbol{\lambda})\in \Delta$ with
$\boldsymbol{\lambda}$ a partition of $N$, 
similar to (\ref{pm}) and (\ref{pmm}), 
we  form sequences
\begin{align*}
& \mbf s(\mbf a^+)= \mbf s(\alpha_0) \cdot \mbf s(\alpha_1) \cdots
\mbf s(\alpha_m)\cdots\\
&\mbf s(\mbf a^-)=\cdots \mbf s(\beta_l) \cdots \mbf s(\beta_{n-1})
\cdot\mbf s(\beta_n)
\end{align*}
where $\alpha_m$ and $\beta_l$ are defined in Section
~\ref{momomial}. 

For each $\mbf a^i$, we denote by $\mbf s(\mbf a^i)$ the sequence
corresponding to $F(O)$ in Proposition ~\ref{inhomogeneous1}.
We set
\begin{align}
\mbf s(\mbf a)= \mbf s(\mbf a^+) \cdot\mbf s(\mbf a^1) \cdots 
\mbf s(\mbf a^s) \cdot\mbf s(N\delta)\cdot \mbf s(\mbf a^-).
\end{align}
The sequence is of finite length since $\mbf a$ is of finite support.
\begin{cor}
\label{singular}
The map $\pi_{\mbf s(\mbf a)}: \tF_{\mbf s(\mbf a)} \to \E_V$ 
(see Section ~\ref{flag}) is a
resolution of singularities in 
$\overline X(\mbf a,\boldsymbol{\lambda})$, 
i.e., $\pi_{\mbf s(\mbf a)}$ is proper,  the image of  $\pi_{\mbf s(\mbf a)}$ is
$\overline X(\mbf a,\boldsymbol{\lambda})$ and the restriction
\[
\pi_{\mbf s(\mbf a)}: 
(\pi_{\mbf s(\mbf a)})^{-1}(X(\mbf a,\boldsymbol{\lambda})) \to 
X(\mbf a,\boldsymbol{\lambda})
\]
is an isomorphism of varieties.
\end{cor}

The proof of Corollary ~\ref{singular} goes exactly the same as the
proof of Theorem 2.2 in ~\cite{Reineke}. The crucial point is Lemma ~\ref{Reineke}. 
We leave the details to the reader.

\end{document}